\pdfoutput=1
\documentclass[toc]{mathprint}
\usepackage{rutarmacros}
\usepackage{macros}
\title[Multifractal Random Substitutions]{Multifractal analysis of measures arising from random substitutions}

\author[Mitchell]
  {Andrew Mitchell}
  {School of Mathematics, University of Birmingham, Edgbaston, B15 2TT, UK}
  {acm925@student.bham.ac.uk}

\author[Rutar]
  {Alex Rutar}
  {Mathematical Institute, University of St Andrews, St Andrews KY16 9SS, Scotland}
  {alex@rutar.org}

\begin{document}
\begin{abstract}
    We study regularity properties of frequency measures arising from random substitutions, which are a generalisation of (deterministic) substitutions where the substituted image of each letter is chosen independently from a fixed finite set.
    In particular, for a natural class of such measures, we derive a closed-form analytic formula for the $L^q$-spectrum and prove that the multifractal formalism holds.
    This provides an interesting new class of measures satisfying the multifractal formalism.
    More generally, we establish results concerning the $L^q$-spectrum of a broad class of frequency measures.
    We introduce a new notion called the \emph{inflation word $L^q$-spectrum} of a random substitution and show that this coincides with the $L^q$-spectrum of the corresponding frequency measure for all $q \geq 0$.
    As an application, we obtain closed-form formulas under separation conditions and recover known results for topological and measure theoretic entropy.
\end{abstract}

\section{Introduction}
A \emph{substitution} is a combinatorial object consisting of a finite alphabet $\mathcal{A}$ along with a set of \emph{transformation rules}.
The theory of substitutions, along with statistical properties of the system under repeated iteration, is a large and actively researched field at the interface of combinatorics and symbolic dynamics.
A thorough introduction to the statistical properties and dynamics of substitutions can be found in \cite{bg2013,que1987}.
Associated with a (deterministic) substitution is a \emph{frequency measure}, which encodes the frequency of subwords under repeated iteration.
Notably, the corresponding subshift supporting this measure has zero topological entropy, and the frequency measure is the unique ergodic measure supported on this subshift.

Random substitutions are a generalisation of (deterministic) substitutions \cite{gl1989} where we apply a transformation rule to each letter randomly and independently chosen from a finite set of possibilities.
Similarly to the deterministic case, subshifts associated with random substitutions support ergodic frequency measures which capture the expected rate of occurrence of subwords under random iteration.
But in contrast to the deterministic case, the corresponding subshift typically has positive topological entropy and supports uncountably many ergodic measures.
Random substitutions include examples exhibiting deterministic behaviour, while also including examples which are subshifts of finite type \cite{grs2019}.
Moreover, there is a large amount of intermediate behaviour: subshifts of random substitutions can simultaneously exhibit long range correlation \cite{bss2018} (an indication of order) and positive topological entropy \cite{goh2020} (an indication of disorder).

As a generalisation of substitutions, random substitutions model \emph{quasicrystals with errors}; namely, crystalline structures with local disorder \cite{al2011,gl1989}.
Random substitutions have also proved a useful tool for modelling fractal percolation \cite{dg1988,dm1990,dv2001}.
Properties of these physical phenomena can be associated with almost sure properties of the frequency measure corresponding to the underlying random substitution \cite{pey1981}.

In this paper, we study the fine scaling properties of frequency measures associated with random substitutions from the perspective of \emph{multifractal analysis}.
This perspective is relevant in a wide variety of contexts, such as the geometry of fractal sets and measures and in dynamical systems, with typical applications to geometric measure theory and number theory.
In our setting, our primary objects of study are the \emph{$L^q$-spectrum}, which is a parametrised family of quantities which capture the inhomogeneous scaling properties of a measure, and the \emph{local dimension}, which capture the exponential growth rate of a measure around a point.
The $L^q$-spectrum and local dimensions are related through a heuristic relationship known as the \emph{multifractal formalism}, first introduced and studied in a physical context in \cite{hjk+1986}.
It is an important and well-studied question to determine settings in which the multifractal formalism holds, and to determine qualitative conditions describing its failure.

Much of the work on multifractal analysis has been done in the setting of local dimensions of self-similar measures (for some examples, see \cite{ap1996,fen2005,fl2009,ln1999,shm2019}) and Birkhoff sums of potentials in dynamical systems with a finite type property (see, for example, \cite{ffw2001,zbl:0933.37020,zbl:0979.28003,pw1999} and the reference therein).
As a notable recent example, in P.~Shmerkin's recent proof of the Furstenberg's intersection conjecture \cite{shm2019}, he computes the $L^q$-spectrum of a large class of dynamically self-similar measures and relates such results to the multifractal analysis of slices of sets.
This information about $L^q$-spectra also implies $L^p$-smoothness properties in the question of absolute continuity of Bernoulli convolutions (see \cite{var2018} for some background on this classic problem).
For more detail on the geometry of measures and multifractal analysis, we refer the reader to the foundational work by L.~Olsen \cite{ols1995} and the classic texts of K.~Falconer \cite{fal1997} and Ya.~Pesin \cite{pes1997}.

Returning to our setting, substitution dynamical systems have characteristic features of (dynamical) self-similarity, but in many cases are far from being ergodic measures on shifts of finite type.
More generally, frequency measures provide a rich family of shift-invariant ergodic measures which exhibit interesting and unique properties in symbolic dynamics in a natural way.
For example, it was proved in \cite[Theorem~4.8]{gmrs2023} that for a certain class of random substitutions, the corresponding subshift supports a frequency measure that is the unique ergodic measure of maximal entropy.
However, this measure is not a Gibbs measure with respect to the zero potential, and therefore the system does not satisfy the common \emph{specification property}, which is a well-known strategy for proving intrinsic ergodicity of symbolic dynamical systems (see \cite{ct2021} and the references therein).
Moreover, there are examples of random substitutions such that the corresponding subshift supports multiple ergodic measures of maximal entropy \cite[Example~6.5]{gs2020}.
More generally, many key properties of frequency measures associated with random substitutions are poorly understood.

In this paper, we derive symbolic expressions for the $L^q$-spectrum of frequency measures associated with random substitutions under certain weak assumptions.
Then under an additional assumption (\textit{recognisability}), we prove a closed-form analytic expression for the $L^q$-spectrum and a variational formula which together imply the multifractal formalism.
Recognisable random substitutions exhibit novel properties not witnessed by other classes of measures for which the multifractal formalism is well-understood: it often happens that the unique frequency measure of maximal entropy is not a Gibbs measure with respect to the zero potential and the corresponding subshift is not sofic.
The techniques and results provide important new perspectives on the geometry and dynamics of the respective measures.

\subsection{Entropy and \texorpdfstring{$L^q$}{Lq}-spectra}
For a Borel probability measure in a compact metric space, the \defn{$L^q$-spectrum} is a well-studied quantity which encodes the scaling properties of the measure, in a weak sub-exponential sense.
Specifically, the $L^q$-spectrum of $\mu$ is given by
\begin{equation*}
    \tau_\mu(q)=\liminf_{r\to 0}\frac{\log \sup\sum_i\mu\bigl(B(x_i,r)\bigr)^q}{\log r}.
\end{equation*}
where the supremum is taken over $2r$-separated subsets $\{x_i\}_i$ of the support of $\mu$.

The $L^q$-spectrum encodes information about the local scaling of the measure $\mu$.
We define the \defn{local dimension} of $\mu$ at $x$ by
\begin{align*}
    \ldim(\mu,x) &=\lim_{r\to 0}\frac{\log\mu\bigl(B(x,r)\bigr)}{\log r}
\end{align*}
when the limit exists.
We then define the \defn{multifractal spectrum} of $\mu$ by
\begin{equation*}
    f_\mu(\alpha) = \dimH\left\{ x \in X : \ldim(\mu,x) = \alpha \right\}\tp
\end{equation*}
In general, the structure of the set of local dimensions can be very complex---for example, the level sets are often dense uncountable subsets of the support of $\mu$.
However, the ``multifractal miracle'' is the phenomenon that, even though the level sets are very complex, the multifractal spectrum is often a concave analytic function of $\alpha$.

In fact, the multifractal spectrum and the $L^q$-spectrum are related through a heuristic relationship called the \emph{multifractal formalism} \cite{hjk+1986}, which speculates that under certain regularity assumptions, the multifractal spectrum is given by the \emph{concave conjugate} of the $L^q$-spectrum, that is the quantity
\begin{equation*}
    \tau_\mu^*(\alpha)=\inf_{q\in\R}(q\alpha-\tau_\mu(q)).
\end{equation*}
Generally speaking, $\tau_\mu^*(\alpha)\geq f_\mu(\alpha)$ \cite{ln1999}: in particular, the slopes of the asymptotes of the $L^q$-spectrum bound the exponential scaling of measures of balls $B(x,r)$ uniformly for all $x\in\supp\mu$.

In our specific setting, where our metric space is the two-sided shift $\mathcal{A}^{\Z}$ and the measure $\mu$ is ergodic, the local dimension is precisely the scaling rate of the \emph{information function} of $\mu$.
In fact, the Shannon--McMillan--Breiman theorem states that the local dimension of the measure (with an appropriate choice of the metric) is almost surely the entropy of the measure.
Thus the $L^q$-spectrum provides \emph{uniform} control over the scaling rate of the information function.
More details about these notions are given in \cref{s:preliminaries}.

\subsection{Random substitutions}
A (deterministic) substitution is a rule which replaces each symbol in a finite or infinite string over an alphabet $\mathcal{A}$ with a finite word over the same alphabet.
Random substitutions generalise this notion by substituting a randomly chosen word (according to a fixed finite distribution) independently for each letter.
We can also think of a random substitution as a (deterministic) set-valued substitution $\vartheta$, together with a choice of probabilities.

For example, given $p \in (0,1)$, the \emph{random Fibonacci substitution} $\vartheta_{p}$ is defined by
\begin{equation*}
    \vartheta_{p} \colon
    \begin{cases}
        a \mapsto
        \begin{cases}
            ab & \text{with probability } p,\\
            ba & \text{with probability } 1-p,
        \end{cases}\\[1.25em]
        b \mapsto a.
    \end{cases}
\end{equation*}
To a given (primitive) random substitution $\vartheta_{\bm{P}}$, one can canonically associate a subshift $X_\vartheta$ of the two-sided shift $\mathcal{A}^{\Z}$ along with an ergodic \emph{frequency measure} $\mu_{\bm{P}}$, which quantifies the relative occurrence of a given word under repeated application of the random substitution.

As highlighted in the introduction, primitive random substitutions give rise to subshifts and measures with a wide variety of properties.
As a result, we impose additional conditions.

Our main assumption, which we call \emph{compatibility} (see \cref{ss:compatibility}), asserts that for each $a\in\mathcal{A}$, the number of occurrences of each $b\in\mathcal{A}$ is identical in every possible substituted image of $a$.
For example, the random Fibonacci substitution is compatible since in all the possible images of $a$, $a$ occurs once and $b$ occurs once.
The key feature of compatibility is that the one can define a deterministic \emph{substitution matrix}, such that the Perron--Frobenius eigenvalue is the asymptotic growth rate of lengths of words under repeated substitution, and the corresponding right eigenvector encodes the asymptotic frequency with which the individual letters appear.
Compatibility is a common assumption: for example, it is assumed in the main results of \cite{bss2018,goh2020,zbl:07720414,rus2020}.

Another standard assumption is \emph{recognisability} (see \cref{ss:recog}), which heuristically states that each element of the subshift has a \emph{unique} preimage in the subshift.
Recognisability precludes the existence of periodic points \cite{rus2020} and is one of the assumptions required to to establish intrinsic ergodicity in \cite{gmrs2023}.
It is also assumed in the main results of \cite{frspreprint,zbl:07720414}.

\subsection{Statement and discussion of main results}
We now give concise statements of the main results in this paper.
We refer the reader to \cref{s:preliminaries} for full statements of the notation and definitions used in this section.

Fix a random substitution $\vartheta_{\bm{P}}$ and let $\lambda$ and $\bm{R}$ denote the Perron--Frobenius eigenvalue and corresponding right eigenvector of the substitution matrix of $\vartheta_{\bm{P}}$, respectively.
Given $q \in \mathbb{R}$ and $k \in \mathbb{N}$, define
\begin{equation*}
    \varphi_{\vartheta_{\bm{P}},k}(q)=\varphi_{k} (q) = - \sum_{a \in \mathcal{A}} R_a \log \left( \sum_{s \in \vartheta^k (a)} \mathbb{P} [\vartheta_{\bm{P}}^k (a) = s]^q \right) \tp
\end{equation*}
We define the \defn{inflation word $L^q$-spectrum of $\vartheta_{\bm{P}}$} by
\begin{equation*}
    T_{\vartheta,\bm{P}}(q) = \liminf_{k \rightarrow \infty} \frac{1}{\lambda^k} \varphi_k(q)\tp
\end{equation*}
We similarly define the upper variant $\lqsymbol{\vartheta,\bm{P}}$ by taking a limit superior in place of the limit inferior.
Throughout, $\mu_{\bm{P}}$ will denote the frequency measure associated with $\vartheta_{\bm{P}}$.
Heuristically, the inflation word spectrum approximates the frequency measure $\mu_{\bm{P}}$ by the probability distribution on the iterated system, with an appropriate normalisation.

Our main general result bounding the $L^q$-spectrum is the following, which states that $T_{\vartheta,\bm{P}}$ and $\tau_{\mu_{\bm{P}}}$ coincide for all $q\geq 0$, and moreover provides bounds on $\tau_{\mu_{\bm{P}}}$ in terms of the functions $\varphi_k$ for all $q\in\R$.
\begin{itheorem}\label{it:lq-inf-equiv}
    Let $\vartheta_{\bm{P}} = (\vartheta, \bm{P})$ be a primitive and compatible random substitution with corresponding frequency measure ${\mu_{\bm{P}}}$.
    Then the limits defining $\tau_{\mu_{\bm{P}}}(q)$ and $T_{\vartheta,\bm{P}}(q)$ exist and coincide for all $q\geq0$.
    Moreover,
    \begin{enumerate}
        \item For all $0\leq q\leq 1$,
            \begin{equation}\label{eq: 0<q<1 bounds}
                \frac{1}{\lambda^k-1} \varphi_k(q) \leq \tau_{\mu_{\bm{P}}}(q) \leq \frac{1}{\lambda^k} \varphi_k(q)
            \end{equation}
            and $(\lambda^{-k}\varphi_k(q))_{k=1}^\infty$ converges monotonically to $\tau_{\mu_{\bm{P}}}(q)$ from above.
        \item For all $q\geq 1$,
            \begin{equation}\label{eq:q>1 bounds}
                \frac{1}{\lambda^k} \varphi_k (q) \leq \tau_{\mu_{\bm{P}}}(q) \leq \frac{1}{\lambda^k-1} \varphi_k(q)
            \end{equation}
            and $(\lambda^{-k}\varphi_k(q))_{k=1}^\infty$ converges monotonically to $\tau_{\mu_{\bm{P}}}(q)$ from below.
        \item For all $q<0$,
            \begin{equation}\label{eq:q<0-lower-bounds}
                \frac{1}{\lambda^k-1} \varphi_k(q)\leq\tau_{\mu_{\bm{P}}}(q)\tp
            \end{equation}
    \end{enumerate}
\end{itheorem}
The notion of compatibility is defined in \cref{ss:compatibility}.
For $q<0$, it is \emph{not} true in general that $\tau_{\mu_{\bm{P}}} (q)$ and $T_{\vartheta,\bm{P}} (q)$ coincide (a counterexample is given in \cref{ex:q<0-non-inflation}): the ``non-uniqueness of cutting points'' allows collisions in the averaging procedure in the construction of the measure (see \cref{l:key-lemma}).
In other words, the corresponding upper bound in \cref{eq:q<0-lower-bounds} does not hold in general.

If $\vartheta_{\bm{P}}$ also satisfies the \defn{disjoint set condition}, or the \defn{identical set condition with identical production probabilities} (see \cref{d:ISC/DSC}), then a closed-form expression can be obtained for $T_{\vartheta,\bm{P}}$ (see \cref{p:inflation-isc-dsc}).
By combining this result with \cref{it:lq-inf-equiv}, we obtain the following corollary.
\begin{icorollary}\label{ic:lq-lim-eq}
    Let $\vartheta_{\bm{P}}$ be a primitive and compatible random substitution with corresponding frequency measure $\mu_{\bm{P}}$.
    Then for all $q\geq 0$:
    \begin{enumerate}
        \item If $\vartheta_{\bm{P}}$ satisfies the disjoint set condition, then
            \begin{equation*}
                \tau_{\mu_{\bm{P}}} (q) = \frac{1}{\lambda - 1} \varphi_1 (q) \tp
            \end{equation*}
        \item If $\vartheta_{\bm{P}}$ satisfies the identical set condition and has identical production probabilities, then
            \begin{equation*}
                \tau_{\mu_{\bm{P}}} (q) = \frac{1}{\lambda} \varphi_1 (q) \tp
            \end{equation*}
    \end{enumerate}
\end{icorollary}
In particular, under the disjoint set condition or identical set condition with identical production probabilities, the $L^q$-spectrum is analytic on $(0,\infty)$.

As a direct application of \cref{it:lq-inf-equiv}, we obtain new proofs of known results on measure-theoretic and topological entropy.
\begin{enumerate}[a]
    \item We obtain the main result of \cite{goh2020} on topological entropy, which states that for subshifts of primitive and compatible random substitutions, the topological entropy can be characterised in terms of the asymptotic growth rate of inflation words.
    \item We also obtain the characterisation of (measure-theoretic) entropy obtained in \cite[Theorem~3.3]{gmrs2023} for frequency measures corresponding to primitive and compatible random substitutions.
        We note that the main results in \cite{gmrs2023} do not require the assumption of compatibility.
\end{enumerate}
This is described in the following corollary.
\begin{icorollary}\label{ic:entropy-recover}
    Let $\vartheta_{\bm{P}} = (\vartheta, \bm{P})$ be a primitive and compatible random substitution with associated subshift $X_{\vartheta}$ and frequency measure $\mu_{\bm{P}}$.
    \begin{enumerate}
        \item The limit
            \begin{equation*}
                \lim_{k \rightarrow \infty} \frac{1}{\lambda^k} \sum_{a \in \mathcal{A}} R_a \log (\# \vartheta^{k} (a))
            \end{equation*}
            exists and is equal to $\htop(X_{\vartheta})$.
        \item The $L^q$-spectrum of $\mu_{\bm{P}}$ is differentiable at $1$.
            Moreover, the limit
            \begin{equation*}
                \lim_{k\to\infty} \frac{1}{\lambda^k} \sum_{a \in \mathcal{A}} R_{a} \sum_{v \in \vartheta^k (a)} - \mathbb{P} [\vartheta_{\bm{P}}^k (a) = v] \log (\mathbb{P} [\vartheta_{\bm{P}}^k (a) = v])
            \end{equation*}
            exists and is equal to $h_{\mu_{\bm{P}}} (X_{\vartheta}) = \dimH\mu_{\bm{P}}=\tau_{\mu_{\bm{P}}}'(1)$.
    \end{enumerate}
\end{icorollary}

We now turn our attention to the multifractal spectrum.
While $\tau_{\mu_{\bm{P}}} (q)$ and $T_{\vartheta,\bm{P}} (q)$ do not coincide in general for $q<0$, if the random substitution that gives rise to the frequency measure $\mu_{\bm{P}}$ is additionally assumed to be recognisable (see \cref{d:recog}), then the limits defining $\tau_{\mu_{\bm{P}}} (q)$ and $T_{\vartheta,\bm{P}} (q)$ both exist and coincide for all $q\in\R$.
Moreover, under recognisability, we prove that the multifractal spectrum is the concave conjugate of the $L^q$-spectrum: in other words, the multifractal formalism holds for any associated frequency measure $\mu_{\bm{P}}$.
In particular, we conclude that $f_{\mu_{\bm{P}}}$ is a concave analytic function.

In fact, in \cref{p:relative-local-dim} we prove a stronger variational formula for the multifractal spectrum.
For each $\alpha\in\R$, we construct measures $\nu$ such that $\dimH \nu\geq\tau^*(\alpha)$ and $\ldim(\mu_{\bm{P}},x)=\alpha$ for $\nu$-a.e. $x\in X_\vartheta$.
In particular, we can take the measures to be frequency measures associated with permissible probabilities for the substitution $\vartheta$.
\begin{itheorem}\label{it:multi-formalism}
    Let $\vartheta_{\bm{P}}$ be a primitive, compatible, and recognisable random substitution with corresponding frequency measure $\mu_{\bm{P}}$.
    Then for all $q\in\R$,
    \begin{equation*}
        \tau_{\mu_{\bm{P}}}(q)=T_{\vartheta,\bm{P}}(q)=\frac{1}{\lambda-1}\varphi_1(q)\tp
    \end{equation*}
    Moreover, $f_{\mu_{\bm{P}}}(\alpha)=\tau_{\mu_{\bm{P}}}^*(\alpha)$ is concave and analytic on its support.
    In fact, for each $\alpha\in\R$ such that $f_{\mu_{\bm{P}}}(\alpha)\geq 0$, there are permissible probabilities $\bm{Q}$ such that $f_{\mu_{\bm{P}}}(\alpha)=\dimH\mu_{\bm{Q}}$ and $\ldim(\mu_{\bm{P}},x)=\alpha$ for $\mu_{\bm{Q}}$-a.e. $x\in X_\vartheta$.
\end{itheorem}

To conclude this section, we observe that our results on $L^q$-spectra also give uniform bounds on the exponential scaling rate of the frequency measures.
The following result is a direct application of \cref{it:lq-inf-equiv} and \cref{it:multi-formalism}, combined with bounds on the possible local dimensions (from \cref{p:meas-bounds}).
\begin{icorollary}\label{ic:ldim-bounds}
    Let $\vartheta_{\bm{P}}=(\vartheta, \bm{P})$ be a primitive, compatible, and recognisable random substitution.
    Then
    \begin{align*}
        \alpha_{\min}&\coloneqq\lim_{q\to\infty}\frac{\tau_{\mu_{\bm{P}}}(q)}{q}=-\sum_{a\in\mathcal{A}}R_a\log\left(\max_{s\in\vartheta(a)}\mathbb{P}[\vartheta_{\bm{P}}(a)=s]\right)\\
        \alpha_{\max}&\coloneqq\lim_{q\to-\infty}\frac{\tau_{\mu_{\bm{P}}}(q)}{q}=-\sum_{a\in\mathcal{A}}R_a\log\left(\min_{s\in\vartheta(a)}\mathbb{P}[\vartheta_{\bm{P}}(a)=s]\right).
    \end{align*}
    and for all $x\in X_{\vartheta}$, $\alpha_{\min}\leq\lldim(\mu_{\bm{P}},x)\leq\uldim(\mu_{\bm{P}},x)\leq\alpha_{\max}$.
    Moreover,
    \begin{equation*}
        \{\ldim(\mu_{\bm{P}},x):x\in X_\vartheta\}=[\alpha_{\min},\alpha_{\max}].
    \end{equation*}
\end{icorollary}
In particular, when the probabilities $\bm{P}$ are chosen so that for each $a\in\mathcal{A}$, $\mathbb{P}[\vartheta_{\bm{P}}(a)=s]=1/\#\vartheta(a)$ for all $s\in\vartheta(a)$, then the $L^q$-spectrum is the line with slope $\htop(X_\vartheta)$ passing through $(1,0)$.
Thus the local dimension of $\mu_{\bm{P}}$ exists at every $x\in X_{\vartheta}$ and is given by the constant value $\alpha_{\min}=\alpha_{\max}$.
This can be rephrased in terms of a weak Gibbs-type property, which says that every sufficiently long legal word (with length depending only on $\epsilon>0$) satisfies
\begin{equation}\label{e:decay-error}
    \exp(-|u|(\htop(X_\vartheta)+\epsilon))\leq \mu_{\bm{P}}([u])\leq \exp(-|u|(\htop(X_\vartheta)-\epsilon))\text{;}
\end{equation}
see, for example, \cite[Lemma~1.4]{shm2019} for a short proof.
In general, the error term $\epsilon$ cannot be dropped by the addition of a constant factor.
Under certain assumptions, one can show that there are infinitely many words with $\mu_{\bm{P}}([u])\approx|u|^{-1}\exp(-|u|(\htop(X_\vartheta))$, as explained in \cite[Lemma~4.12]{gmrs2023}.
These assumptions are satisfied, for example, in \cref{ex:intrinsic}.

Of course, similar one-sided results hold for $q\geq 0$ only under the assumption of compatibility, by iterating the formula for $\varphi_k$ and taking an appropriate maximum at each level.
In fact, since $\tau_{\mu_{\bm{P}}}(q)$ is differentiable at $1$, with derivative giving the entropy, and since $\htop(X_\vartheta)=\tau_{\mu_{\bm{P}}}(0)$, it follows that $\mu_{\bm{P}}$ is a measure of maximal entropy if and only if $\tau_{\mu_{\bm{P}}}'(q)$ exists and is constant on the interval $(0,1)$.

\subsection{Further work}
We conclude the introduction with a list of comments and potentially interesting questions.
\begin{enumerate}
    \item What is the $L^q$-spectrum for a compatible substitution when $q<0$?
        We do not know this even for the random substitution given in \cref{ex:q<0-failure-full-shift}, which satisfies the identical set condition with identical production probabilities.
        Obtaining results for $q<0$ is substantially more challenging, since the sum in $\tau_{\mu_{\bm{P}}}(q)$ depends on the measure of cylinders with very small (but non-zero) measure.
        For example, in the self-similar case, without the presence of strong separation assumptions, little is known (in contrast to the $q\geq 0$ case).
    \item What happens without compatibility?
        Do the formulas in \cref{it:lq-inf-equiv} hold in general?
        In \cite{gmrs2023}, it suffices to use an almost sure version of \cref{l:comp-letter-freqs}.
        However, since the $L^q$-spectrum is sensitive to scaling at individual points as $q$ tends to $\pm \infty$, such an almost sure result in our case is insufficient.
    \item Outside the disjoint set condition and the identical set condition, what can be said about differentiability of the $L^q$-spectrum?
        For $q\geq 0$, we give the $L^q$-spectrum as a uniform limit of analytic functions: however, aside from the exceptional point $q=1$ where we can say more, this is not enough to give information about differentiability.
    \item Can the assumption of recognisability in \cref{it:multi-formalism} be relaxed to a weaker condition, such as the disjoint set condition (see \cref{d:ISC/DSC})?
    \item Can the error term in \cref{e:decay-error} be determined precisely, up to a constant?
        The approximate Gibbs-type bounds discussed following \cref{ic:ldim-bounds} are closely related to the bounds used in the proof of intrinsic ergodicity given in \cite[Theorem~4.8]{gmrs2023}.
        It could be worth exploring the relationship between intrinsic ergodicity and Gibbs-type properties given by the $L^q$-spectrum.
\end{enumerate}

\section{Preliminaries}\label{s:preliminaries}
In this section we introduce the key notation and definitions that we will use throughout the paper.
After introducing some basic notation, in \cref{ss:dynamics} we introduce symbolic dynamics on the two-sided shift, as well as our notions of entropy and dimension.
In \cref{ss:lq-spectra} we present the key definitions and basic results from multifractal analysis that we work with throughout, including the definitions of the $L^q$-spectrum and local dimensions of a measure.
Then, in the following sections we provide an introduction to random substitutions and their associated dynamical systems.
In \cref{ss:random-subst} we give the definition of a random substitution via its action on words, and define the subshift associated to a random substitution.
Then, in \cref{ss:primitivity} and \cref{ss:compatibility}, we define what it means for a random substitution to be \defn{primitive} and \defn{compatible} and present the key properties of such random substitutions.
In \cref{ss:freq-measures}, we give the definition of the frequency measure associated to a random substitution and state a key result used in the proof of our main results which relates the measures of cylinder sets via the action of the random substitution.
Finally, in \cref{ss:recog}, we define what it means for a substitution to satisfy the disjoint or identical set condition, and introduce recognisable random substitutions.

\subsection{Symbolic notation}
Throughout, we use the following symbolic notation, which is essentially the same as the notation used in \cite{bg2013,lm1995}.

For a set $B$, we let $\# B$ be the cardinality of $B$ and let $\mathcal{F}(B)$ be the set of non-empty finite subsets of $B$.

We fix an \defn{alphabet} $\mathcal{A} = \{ a_{1}, \ldots, a_{d} \}$, for some $d \in \mathbb{N}$, which is a finite set of \defn{letters} $a_{i}$, and equip it with the discrete topology.
Then a \defn{word} $u$ with letters in $\mathcal{A}$ is a finite concatenation of letters, namely $u = a_{i_{1}} \cdots a_{i_{n}}$ for some $n \in \mathbb{N}$.
We write $\lvert u \rvert = n$ for the length of the word $u$, and for $m \in \mathbb{N}$, we let $\mathcal{A}^{m}$ denote the set of all words of length $m$ with letters in $\mathcal{A}$.

We set $\mathcal{A}^{+} = \bigcup_{m \in \mathbb{N}} \mathcal{A}^{m}$ and let
\begin{equation*}
    \mathcal{A}^{\mathbb{Z}} = \{ (a_{i_n})_{n\in\Z} : a_{i_n} \in \mathcal{A} \text{ for all }n \in \mathbb{Z} \}
\end{equation*}
denote the set of all bi-infinite sequences with elements in $\mathcal{A}$, and endow $\mathcal{A}^{\mathbb{Z}}$ with the product topology.
We also fix a metric on $\mathcal{A}^{\Z}$ as follows.
Given points $x=(x_n)_{n\in\Z}$ and $y=(y_n)_{n\in\Z}$, let $N(x,y)=\sup\{n\in\Z:x_j=y_j\text{ for all }|j|\leq n\}$ and let $d(x,y)=e^{-2N(x,y)-1}$.
The space $X$ is compact with topology generated by the metric.

We will frequently write sequences $(x_n)_{n\in\Z}\in\mathcal{A}^{\Z}$ as $\cdots x_{-1}x_0x_1\cdots$, with the corresponding notation for finite sequences.
If $i$ and $j \in \mathbb{Z}$ with $i \leq j$, and $x = \cdots x_{-1} x_{0} x_{1} \cdots \in \mathcal{A}^{\mathbb{Z}}$, then we let $x_{[i,j]} = x_i x_{i+1} \cdots x_{j}$.
We use the same notation if $v \in \mathcal{A}^{+}$ and $1 \leq i \leq j \leq |v|$.
For $u$ and $v \in \mathcal{A}^{+}$ (or $v \in \mathcal{A}^{\Z}$), we write $u \triangleleft v$ if $u$ is a subword of $v$, namely if there exist $i$ and $j \in \mathbb{Z}$ with $i \leq j$ so that  $u = v_{[i, j]}$.
For $u$ and $v \in \mathcal{A}^{+}$, we set $\lvert v \rvert_u$ to be the number of (possibly overlapping) occurrences of $u$ as a subword of $v$.
If $u = a_{i_1} \cdots a_{i_n}$ and $v = a_{j_1} \cdots a_{j_m} \in \mathcal{A}^{+}$, for some $n$ and $m \in \mathbb{N}$, we write $u v$ for the \defn{concatenation} of $u$ and $v$.
The \defn{abelianisation} of a word $u \in \mathcal{A}^{+}$ is the vector $\Phi (u) \in \N_0^{\#\mathcal{A}}$, defined by $\Phi (u)_{a} = \lvert u \rvert_{a}$ for all $a\in\mathcal{A}$.

\subsection{Dynamics, entropy and dimension}\label{ss:dynamics}
We equip the space $\mathcal{A}^{\Z}$ with invertible dynamics from the \defn{left-shift map} $S\colon\mathcal{A}^{\Z}\to\mathcal{A}^{\Z}$.
Throughout, we work with a \defn{subshift} $X\subset\mathcal{A}^{\Z}$, which is compact and \defn{shift-invariant}, that is $S^{-1}(X)=X$.
Then $\mu$ will denote an ergodic and $S$-invariant Borel probability measure with support contained in $X$.

The metric structure on $\mathcal{A}^{\Z}$ enables us to define the \defn{Hausdorff dimension} of Borel subsets of $X$.
Using this, we define the \defn{Hausdorff dimension of $\mu$} to be the quantity
\begin{equation*}
    \dimH\mu=\inf\{\dimH E:\mu(E)>0\}
\end{equation*}
where the infimum is taken over Borel sets $E$.
Here, we use Hausdorff dimension as inherited from the underlying metric; though it would also be appropriate to use Bowen's generalisation of topological entropy to non-compact sets \cite{bow1973}.
We also define the \defn{lower local dimension} of $\mu$ at $x$ by
\begin{align*}
    \lldim(\mu,x) &=\liminf_{r\to 0}\frac{\log\mu\bigl(B(x,r)\bigr)}{\log r}
\end{align*}
We define the \defn{upper local dimension} $\uldim(\mu,x)$ analogously using the limit superior in place of the limit inferior, and when the limits coincide, we refer to the shared quantity as the \defn{local dimension} and denote it by $\ldim(\mu,x)$.

Local dimensions and Hausdorff dimension are closely related: the same proof as given, for instance, in \cite[Proposition~10.1]{fal1997} implies that
\begin{equation}\label{e:Hdim-ldim}
    \dimH\mu=\sup\{s:\lldim(\mu,x)\geq s\text{ for $\mu$-a.e. }x\}.
\end{equation}

Now fix a partition $\xi$ so that with $\xi_k=\bigvee_{i=-k}^k S^{-i}(\xi)$, $\{\xi_k\}_{k=1}^\infty$ generates the Borel $\sigma$-algebra on $X$.
We recall that the \defn{entropy} of $\mu$ with respect to $S$ is given by
\begin{equation*}
    h_\mu(X)=\lim_{k\to\infty}\frac{1}{2k+1}\sum_{A\in\xi_k}-\mu(A)\log\bigl(\mu(A)).
\end{equation*}
where, by the classical Kolmogorov--Sinaĭ theorem, the quantity does not depend on the choice of partition.

Now given $x\in X$, let $\xi_k(x)$ denote the unique element in the partition $\xi_k$ containing $x$.
Then the Shannon--McMillan--Breiman theorem states that the entropy of $\mu$ is almost surely the \defn{information rate} of $\mu$, that is for $\mu$-a.e. $x\in X$,
\begin{equation*}
    \lim_{k\to\infty}\frac{-\log\mu\bigl(\xi_k(x)\bigr)}{2k+1}=h_\mu(X).
\end{equation*}
We refer the reader to \cite[Theorem~3.2.7]{zbl:0896.28006} for a proof and more background on this topic.

Now suppose $\xi=\{E_a\}_{a\in\mathcal{A}}$ is the partition of $X$ where $E_a=\{(x_n)_{n\in\Z}\in X:x_0=a\}$.
For the remainder of this paper, $\xi$ will always denote this partition.
Then given $x=(x_n)_{n\in\Z}\in X$,
\begin{equation*}
    \xi_k(x)=\{y\in X:x_j=y_j\text{ for all }|j|\leq k\}=B(x,e^{-(2k+1)})
\end{equation*}
and therefore
\begin{equation*}
    \ldim(\mu,x) = \lim_{k\to\infty}\frac{-\log\mu\bigl(\xi_k(x)\bigr)}{2k+1}
\end{equation*}
where both limits exist if either limit exists.
Since the limit on the right is $\mu$ almost surely $h_\mu(X)$, it follows from \cref{e:Hdim-ldim} that $\dimH\mu=h_\mu(X)$.

Finally, the \defn{topological entropy} of $X$ is given by
\begin{equation*}
    \htop(X)=\lim_{k\to\infty}\frac{-\log\#\{E\in\xi_k:E\cap X\neq\varnothing\}}{2k+1}.
\end{equation*}
Of course, $\htop(X)=\dimB X$, the box counting dimension of $X$.

\subsection{\texorpdfstring{$L^q$}{Lq}-spectra and smoothness}\label{ss:lq-spectra}
Given $q\in\R$, we define
\begin{equation*}
    \overline{S}_{\mu,r}(q)=\sup_{\{x_i\}_i\in\mathcal{P}(r)}\sum_{i}\mu\bigl(B(x_i,r)\bigr)^q
\end{equation*}
where $\mathcal{P}(r)$ is the set of discrete $2r$-separated subsets of $X$, that is $\mathcal{P}(r)=\{\{x_i\}_{i}:x_i\in X,d(x_i,x_j)\geq 2r\text{ for }i\neq j\}$.
We then define the \defn{$L^q$-spectrum of $\mu$} to be the function
\begin{equation*}
    \tau_\mu(q)=\liminf_{q\to 0}\frac{\log \overline{S}_{\mu,r}(q)}{\log r}.
\end{equation*}
For convenience, we also denote the upper variant $\overline{\tau}_\mu(q)$ by taking a limit superior in place of the limit inferior.
It is a standard consequence of Hölder's inequality that $\tau_\mu(q)$ is a concave increasing function of $q$ (note that this need not hold for $\overline{\tau}_\mu(q)$).

Of course, the preceding definitions hold more generally in an arbitrary metric space, but in our particular setting we can rephrase the $L^q$-spectrum in terms of more familiar sums over cylinders.
Recall that $\xi$ denotes the partition of $X_\vartheta$ into cylinders at $0$ corresponding to the letters in $\mathcal{A}$.
Then set
\begin{equation*}
    S_{\mu,n}(q)=\sum_{E\in \xi_k}\mu(E)^q.
\end{equation*}
Since distinct elements in the partition $\xi_k$ are $e^{-(2k+1)}$-separated,
\begin{equation*}
    S_{\mu,n}(q)=\overline{S}_{\mu,e^{-(2n+1)}}(q).
\end{equation*}
It follows immediately that
\begin{equation*}
    \tau_\mu(q)=\liminf_{n\to\infty}\frac{-\log S_{\mu,n}(q)}{2n+1}
\end{equation*}
with the analogous result for $\overline{\tau}_\mu(q)$.
In particular, $\tau_\mu(0)=\htop(X)$ assuming $\mu$ is fully supported on $X$.

Finally, by shift invariance, we can characterise the subshift $X$ in terms of a \defn{language} on $X$.
Given $n\in\N$, we set
\begin{equation*}
    \mathcal{L}^n(X)=\{w\in\mathcal{A}^n:w\triangleleft x\text{ for some }x\in X\}.
\end{equation*}
Given $w\in\mathcal{L}^n(X)$, we let $[w]=\{(x_n)_{n\in\Z}\in X:x_i=w_i\text{ for all }1\leq i\leq n\}$.
Of course, by shift invariance, there is a measure-preserving bijection between $\mathcal{L}^{2n+1}(X)$ and $X_n$, so it follows again that
\begin{equation*}
    \tau_\mu(q)=\liminf_{n\to\infty} -\frac{1}{n} \log\sum_{u\in\mathcal{L}^n(X)}\mu([u])^q.
\end{equation*}
We will primarily use this characterisation throughout the paper.

We first list some basic properties of the $L^q$-spectrum of the measure $\mu$.
Here, (a) is a direct consequence of Hölder's inequality, (b) is standard (see, for example, \cite[Lemma~1.4]{shm2019}) and (c) was proved in \cite[Theorem~1.4]{flr2002}.
\begin{lemma}\label{l:lq-core}
    Let $\mu$ be a shift-invariant measure on $X$.
    \begin{enumerate}[a]
        \item The $L^q$-spectrum $\tau_\mu(q)$ is continuous, increasing and concave on $\R$.
        \item Let $\alpha_{\min}=\lim_{q\to\infty}\tau_\mu(q)/q$ and $\alpha_{\max}=\lim_{q\to -\infty}\tau_\mu(q)/q$.
            Then for every $s<\alpha_{\min}\leq\alpha_{\max}<t$, all $n$ sufficiently large and $u\in\mathcal{L}^n$, $e^{-tn}\leq \mu([u])\leq e^{-s n}$.
            In particular, the local dimensions satisfy
            \begin{equation*}
                \alpha_{\min}\leq\inf_{x\in X}\lldim(\mu,x)\leq\sup_{x\in X}\uldim(\mu,x)\leq\alpha_{\max}.
            \end{equation*}
        \item The left and right derivatives of $\tau_\mu$ at $q=1$ bound the Hausdorff dimension of $\mu$, that is $\tau_\mu^+(1)\leq \dimH\mu\leq\tau_\mu^-(1)$.
    \end{enumerate}
\end{lemma}
In fact, (a) gives intuition for why the $L^q$-spectrum encodes smoothness: rather than obtain almost sure information on local dimensions, the $L^q$-spectrum contains \emph{uniform} information about local dimensions.

Finally, we prove a simple result concerning the $L^q$-spectrum which will be useful later in the paper.
\begin{lemma}\label{l:lq-subseq}
    Let $\zeta>1$ be arbitrary.
    Then
    \begin{equation}\label{e:lq-lower-sub}
        \tau_{\mu}(q)=\frac{1}{\zeta}\liminf_{n\to\infty} - \frac{1}{n} \log \left( \sum_{u \in \mathcal{L}^{\lfloor \zeta n\rfloor}(X)} \mu ([u])^q \right)
    \end{equation}
    and
    \begin{equation}\label{e:lq-upper-sub}
        \lqupper{\mu}(q)=\frac{1}{\zeta}\limsup_{n\to\infty} - \frac{1}{n} \log \left( \sum_{u \in \mathcal{L}^{\lfloor \zeta n\rfloor} (X)} \mu ([u])^q \right)\tp
    \end{equation}
\end{lemma}
\begin{proof}
    Of course, it always holds that
    \begin{align*}
        \tau_\mu(q)&\leq\frac{1}{\zeta}\liminf_{n\to\infty} - \frac{1}{n} \log \left( \sum_{u \in \mathcal{L}^{\lfloor \zeta n\rfloor} (X)} \mu ([u])^q \right)\\
        \lqupper{\mu}(q)&\geq\frac{1}{\zeta}\limsup_{n\to\infty} - \frac{1}{n} \log \left( \sum_{u \in \mathcal{L}^{\lfloor \zeta n\rfloor} (X)} \mu ([u])^q \right)\tp
    \end{align*}
    First, let $q<0$ and let $n\in\N$ be arbitrary.
    Let $k_n$ be minimal so that $\lfloor \zeta k_n\rfloor\geq n$.
    Observe that there is some $M\in\N$ (independent of $n$) so that $\lfloor\zeta k_n\rfloor\leq n+M$: it follows that $\lim_{n\to\infty} n/k_n=\zeta$.
    Then if $v\in\mathcal{L}^{\lfloor\zeta k_n\rfloor}(X)$ is arbitrary, $[v]\subset[u]$ for some $u\in\mathcal{L}^n(X)$ and $\mu([v])^q\geq\mu([u])^q$.
    Thus
    \begin{equation*}
        S_{\lfloor\zeta k_n\rfloor,\mu}(q)\geq S_{n,\mu}(q)\tp
    \end{equation*}
    which gives \cref{e:lq-lower-sub} for $q<0$ since $\lim n/k_n=\zeta$.

    Similarly, for $q\geq 0$, since there are at most $(\#\mathcal{A})^M$ words $v\in\mathcal{L}^{\lfloor\zeta k_n\rfloor}(X)$ with $[v]\subset[u]$, pigeonholing, for each $u\in\mathcal{L}^n(X)$ there is some $v(u)\in\mathcal{L}^{\lfloor\zeta k_n\rfloor}(X)$ such that $\mu([v(u)])^q\geq (\#\mathcal{A})^{-qM}\mu([u])^q$.
    Thus
    \begin{equation*}
        S_{\lfloor\zeta k_n\rfloor,\mu}(q)\geq (\#\mathcal{A})^{-qM} S_{n,\mu}(q)\tp
    \end{equation*}
    This gives \cref{e:lq-lower-sub} for $q\geq 0$.

    The arguments for \cref{e:lq-upper-sub} follow analogously by choosing $k_n$ maximal so that $\lfloor\zeta k_n\rfloor\leq n$.
\end{proof}

\subsection{Multifractal spectrum and multifractal formalism}
The $L^q$-spectrum of a measure is related to the \defn{(fine) multifractal spectrum}.
Let $\mu$ be a shift-invariant measure on a subshift $X$.
We recall that the local dimension of $\mu$ at $x\in X$ is given by
\begin{align*}
    \ldim(\mu,x) &= \lim_{n \rightarrow \infty} - \frac{1}{2n+1} \log \mu ([x_{[-n,n]}])
\end{align*}
when the limit exists.
Given $\alpha \in \mathbb{R}$, set
\begin{equation*}
    F_\mu(\alpha) = \left\{ x \in X : \ldim(\mu,x) = \alpha \right\} \tp
\end{equation*}
We then define the \defn{multifractal spectrum} of $\mu$ by
\begin{equation*}
    f_{\mu}(\alpha) = \dimH F_\mu(\alpha)
\end{equation*}
using the convention that $\dimH\varnothing=-\infty$.

The multifractal spectrum is related to the $L^q$-spectrum by the following result.
Let $g\colon\R\to\R\cup\{-\infty\}$ be a concave function.
For $x\in\R$, we let $g^+(x)$ (resp. $g^-(x)$) denote the right (resp. left) derivative of $g$ at $x$.
Such limits necessarily exist by concavity.
We denote the \defn{subdifferential} of $g$ at $x$ by $\partial g(x)=[g^+(x),g^-(x)]$.
We then recall that the \defn{concave conjugate} of $g$ is given by
\begin{equation*}
    g^*(\alpha)=\inf_{q\in\R}\{q\alpha-g(q)\}.
\end{equation*}
Note that $g^*$ is always concave since it is the infimum of a family of affine functions.
For more detail concerning the theory of concave functions, we refer the reader to \cite{roc1970}.

Now, we say that $\mu$ satisfies the \defn{multifractal formalism} when $f_{\mu}=\tau_\mu^*$.
In general, the multifractal formalism need not hold, but it is well-known that the concave conjugate of the $L^q$-spectrum is an upper bound for the multifractal spectrum.
For the convenience of the reader, we provide a short self-contained proof, which follows the main ideas of \cite[Theorem~4.1]{ln1999}.
\begin{proposition}\label{p:multi-upper}
    Let $\mu$ be a shift-invariant measure on a subshift $X$.
    Then $f_\mu(\alpha)\leq\tau^*(\alpha)$ for all $\alpha\in\R$.
\end{proposition}
\begin{proof}
    Recall that $\xi_n$ denotes the partition of $X$ into cylinders corresponding to words of length $2n+1$, each of which has diameter precisely $e^{-(2n+1)}$.
    For $\alpha\in\R$, $n\in\N$ and $\epsilon>0$, let
    \begin{equation*}
        \mathcal{M}_{n,\epsilon}(\alpha)=\Bigl\{I\in\xi_n:e^{-(2n+1)(\alpha+\epsilon)}\leq\mu(I)\leq e^{-(2n+1)(\alpha-\epsilon)}\Bigr\}\tp
    \end{equation*}
    In other words, $\mathcal{M}_{n,\epsilon}(\alpha)$ is an $\epsilon$-approximation of $F_\mu(\alpha)$ at level $n$.
    Our strategy is to control the size of the sets $\mathcal{M}_{n,\epsilon}(\alpha)$ in terms of the $L^q$-spectrum of $\mu$, and then use these sets to build a good cover of $F_\mu(\alpha)$.
    Let $q\in\partial\tau^*(\alpha)$: we prove this in the case that $q\geq 0$; the case $q<0$ is analogous.

    First,
    \begin{equation}\label{e:sn-upper}
        S_{2n+1,\mu}(q)=\sum_{I\in\xi_n}\mu(I)^q\geq\sum_{u\in\mathcal{M}_{n,\epsilon}(\alpha)}\mu(I)^q\geq e^{-(2n+1)(\alpha+\epsilon)q}\#\mathcal{M}_{n,\epsilon}(\alpha)\tp
    \end{equation}
    Since $\tau_\mu(q)=\liminf_{n\to\infty}(\log S_{2n+1,\mu}(q))/(-2n-1)$ by \cref{l:lq-subseq}, there is some $N_\epsilon\in\N$ so that for all $n\geq N_\epsilon$, $S_{2n+1,\mu}(q)\leq e^{-(2n+1)(\tau_\mu(q)-\epsilon)}$.
    Combining this with \cref{e:sn-upper},
    \begin{equation}\label{e:M-cover-bound}
        \#\mathcal{M}_{n,\epsilon}(\alpha)\leq e^{-(2n+1)(\tau(q)-\epsilon)}\cdot e^{(2n+1)(\alpha+\epsilon)q}= e^{(2n+1)(\tau^*(\alpha)+(q+1)\epsilon)}
    \end{equation}
    for all $n\geq N_\epsilon$ where we have used the fact that $q\in\partial\tau^*(\alpha)$.

    Now for each $x\in F_\mu(\alpha)$, we can find some $n_x\in\N$ so that for all $n\geq n_x$, $\mu(\xi_n(x))\geq e^{-(2n+1)(\alpha+\epsilon)}$.
    In particular,
    \begin{equation*}
        \mathcal{G}_\epsilon\coloneqq\bigcup_{n=N_\epsilon}^\infty\mathcal{M}_{n,\epsilon}(\alpha)
    \end{equation*}
    is a Vitali cover for $F_\mu(\alpha)$.

    Now suppose $\{I_j\}_{j=1}^\infty$ is any disjoint subcollection of $\mathcal{G}_\epsilon$: then with $s=\tau^*(\alpha)+2\epsilon(1+q)$,
    \begin{align*}
        \sum_{j=1}^\infty(\diam I_j)^s&\leq\sum_{n=N_\epsilon}^\infty\sum_{I\in\mathcal{M}_{n,\epsilon}(\alpha)}(\diam I)^s\leq\sum_{n=N_\epsilon}^\infty e^{-(2n+1)s}\#\mathcal{M}_{n,\epsilon}(\alpha)\\
                                        &\leq\sum_{n=N_\epsilon}^\infty e^{-(2n+1)s}e^{(2n+1)(\tau^*(\alpha)+(q+1)\epsilon)}\\
                                        &=\sum_{n=N_\epsilon}^\infty (e^{-(1+q)\epsilon})^{2n+1}<\infty
    \end{align*}
    by \cref{e:M-cover-bound}.
    Thus by the Vitali covering theorem for Hausdorff measure (\cite[Theorem~1.10]{zbl:0587.28004} holds in the shift setting with the same proof), there is a cover $\{E_i\}_{i=1}^\infty$ for $F_\mu(\alpha)$ such that
    \begin{equation*}
        \mathcal{H}^s(F_\mu(\alpha))\leq \sum_{i=1}^\infty(\diam E_i)^s<\infty
    \end{equation*}
    and thus $\dimH F_\mu(\alpha)\leq\tau^*(\alpha)+2\epsilon(1+q)$.
    But $\epsilon>0$ was arbitrary, so the desired result follows.
\end{proof}

\subsection{Random substitutions}\label{ss:random-subst}
We now introduce our primary objects of interest: \defn{random substitutions}, and their associated \defn{frequency measures}.
In a similar manner to \cite{goh2020,gmrs2023}, we define a random substitution by the data required to determine its action on letters.
We then extend this to a random map on words.
\begin{definition}
    Let $\mathcal{A} = \{ a_{1}, \ldots, a_{d} \}$ be a finite alphabet.
    A \defn{random substitution} $\vartheta_{\bm{P}} = (\vartheta, \bm{P})$ is a \defn{set-valued substitution} $\vartheta \colon \mathcal{A} \rightarrow \mathcal{F}(\mathcal{A}^{+})$ together with a set of non-degenerate probability vectors $\bm{P}=(\bm{p}_i)_{i=1}^d$ where
    \begin{equation*}
        \bm{p}_{i} = ( p_{i, 1}, \ldots, p_{i, r_i} )\qquad\text{with}\qquad r_i = \# \vartheta(a_i); \, \bm{p}_i \in (0,1]^{r_i}; \, \sum_{j=1}^{r_i} p_{i,j} = 1.
    \end{equation*}
    such that
    \begin{align*}
        \vartheta_{\bm{P}} \colon a_i \mapsto
        \begin{cases}
            s^{(i,1)} & \text{with probability } p_{i, 1},\\
            \hfill \vdots \hfill & \hfill \vdots\hfill\\
            s^{(i,r_i)} & \text{with probability } p_{i, r_i},
        \end{cases}
    \end{align*}
    for every $1 \leq i \leq d$, where $\vartheta(a_i) = \{ s^{(i,j)} \}_{1\leq j \leq r_i}$.

    We call each $s^{(i,j)}$ a \defn{realisation} of $\vartheta_{\bm{P}}(a_i)$.
    If $r_i = 1$ for all $i \in \{ 1, \ldots, d \}$, then we call $\vartheta_{\bm{P}}$ \defn{deterministic}.
\end{definition}
\begin{example}[Random Fibonacci]\label{ex:random-fib}
    Let $\mathcal{A} = \{ a, b\}$, and let $p \in (0,1)$.
    The \defn{random Fibonacci substitution} $\vartheta_{\bm{P}} = (\vartheta, \bm{P})$ is the random substitution given by
    \begin{align*}
        \vartheta_{\bm P} \colon
        \begin{cases}
            a \mapsto
            \begin{cases}
                ab & \text{with probability } p,\\
                ba & \text{with probability } 1-p,
            \end{cases}\\[1.25em]
            b \mapsto a
        \end{cases}
    \end{align*}
    with defining data $r_{a} = 2$, $r_{b} = 1$, $s^{(a, 1)} = ab$, $s^{(a, 2)} = ba$, $s^{(b, 1)} = a$, $\bm{P} = \{ \bm{p}_{a} = (p, 1-p), \bm{p}_{b} = (1) \}$ and corresponding set-valued substitution $\vartheta \colon a \mapsto \{ab,ba\}, b \mapsto \{a\}$.
\end{example}

In the following we describe how a random substitution $\vartheta_{\bm{P}}$ determines a (countable state) Markov matrix $Q$, indexed by $\mathcal{A}^{+} \times \mathcal{A}^{+}$.
We interpret the entry $Q_{u,v}$ as the probability of mapping a word $u$ to a word $v$ under the random substitution.
Formally, $Q_{a_i, s^{(i,j)}} = p_{i,j}$ for all $j \in \{1,\ldots, r_i\}$ and $Q_{a_i,v} =0$ if $v \notin \vartheta(a_i)$.
We extend the action of $\vartheta_{\bm{P}}$ to finite words by mapping each letter \emph{independently} to one of its realisations.
More precisely, given $n \in \mathbb{N}$, $u = a_{i_1} \cdots a_{i_n} \in \mathcal{A}^{n}$ and $v \in \mathcal{A}^{+}$ with $|v| \geq n$, we let
\begin{align*}
    \mathcal{D}_n(v) = \{ (v^{(1)},\ldots, v^{(n)}) \in (\mathcal{A}^{+})^{n} : v^{(1)} \cdots v^{(n)} = v \}
\end{align*}
denote the set of all decompositions of $v$ into $n$ individual words and set
\begin{align*}
    Q_{u,v} = \sum_{(v^{(1)},\ldots,v^{(n)}) \in \mathcal{D}_n(v)} \prod_{j = 1}^{n} Q_{a_{i_j},v^{(j)}}.
\end{align*}
In words, $\vartheta_{\bm{P}}(u) = v$ with probability $Q_{u,v}$.

For $u \in \mathcal{A}^{+}$, let $(\vartheta_{\bm{P}}^{n}(u))_{n \in \mathbb{N}}$ be a stationary Markov chain on some probability space $(\Omega_u, \mathcal{F}_u, \mathbb{P}_u)$, with Markov matrix given by $Q$; that is,
\begin{align*}
    \mathbb{P}_u [\vartheta_{\bm{P}}^{n+1}(u) = w \mid \vartheta_{\bm{P}}^{n}(u) = v] = \mathbb{P}_v [\vartheta_{\bm{P}}(v) = w] = Q_{v,w}
\end{align*}
for all $v$ and $w \in \mathcal{A}^{+}$, and $n \in \mathbb{N}$.
In particular,
\begin{align*}
    \mathbb{P}_u [\vartheta_{\bm{P}}^{n}(u) = v] = (Q^{n})_{u,v}
\end{align*}
for all $u$ and $v \in \mathcal{A}^{+}$, and $n \in \mathbb{N}$.
We often write $\mathbb{P}$ for $\mathbb{P}_u$ if the initial word is understood.
In this case, we also write $\mathbb{E}$ for the expectation with respect to $\mathbb{P}$.
As before, we call $v$ a \defn{realisation} of $\vartheta^{n}_{\bm{P}}(u)$ if $(Q^{n})_{u,v} > 0$ and set
\begin{align*}
    \vartheta^{n}(u) = \{ v \in \mathcal{A}^{+} : (Q^{n})_{u,v} > 0\}
\end{align*}
to be the set of all realisations of $\vartheta_{\bm{P}}^{n}(u)$.
Conversely, we may regard $\vartheta^{n}_{\bm{P}}(u)$ as the set $\vartheta^{n}(u)$ endowed with the additional structure of a probability vector.
If $u = a \in \mathcal{A}$ is a letter, we call a word $v \in \vartheta^{k}(a)$ a \defn{level-$k$ inflation word}, or \defn{exact inflation word}.

To a given random substitution $\vartheta_{\bm{P}} = (\vartheta, \bm{P})$ one can associate a subshift.
First, we say that a word $u \in \mathcal{A}^{+}$ is \defn{($\vartheta$-)legal} if there exists an $a_i \in \mathcal{A}$ and $k \in \mathbb{N}$ such that $u$ appears as a subword of some word in $\vartheta^{k} (a_i)$.
We define the \defn{language} of $\vartheta$ by $\mathcal{L}_{\vartheta} = \{ u \in \mathcal{A}^{+} : u \text{ is $\vartheta$-legal} \}$ and, for $w \in \mathcal{A}^{+} \cup \mathcal{A}^{\mathbb{Z}}$, we let $\mathcal{L} (w) = \{ u \in \mathcal{A}^{+} : u \triangleleft w \}$ denote the language of $w$.
\begin{definition}
    The \defn{random substitution subshift} of a random substitution $\vartheta_{\bm{P}} = (\vartheta, \bm{P})$ is the system $(X_{\vartheta}, S)$, where $X_{\vartheta} = \{ w \in \mathcal{A}^{\mathbb{Z}} : \mathcal{L} (w) \subseteq \mathcal{L}_{\vartheta} \}$ and $S$ denotes the (left) shift map, defined by $S (w)_{i} = w_{i+1}$ for each $w \in X_{\vartheta}$.
\end{definition}

Under very mild assumptions, the space $X_{\vartheta}$ is non-empty \cite{gs2020}.
This holds, for example, if the generating random substitution is primitive (we give a definition in \cref{ss:primitivity}).
We endow $X_{\vartheta}$ with the subspace topology inherited from $\mathcal{A}^{\mathbb{Z}}$, and since $X_{\vartheta}$ is defined in terms of a language, it is a compact $S$-invariant subspace of $\mathcal{A}^{\mathbb{Z}}$.
Hence, $X_{\vartheta}$ is a subshift.
For $n \in \mathbb{N}$, we write $\mathcal{L}_{\vartheta}^{n} = \mathcal{L}_\vartheta \cap \mathcal{A}^{n}$ and $\mathcal{L}^{n} (w) = \mathcal{L}(w) \cap \mathcal{A}^{n}$ to denote the subsets of $\mathcal{L}_{\vartheta}$ and $\mathcal{L} (w)$, respectively, consisting of words of length $n$.
The set-valued function $\vartheta$ naturally extends to $X_{\vartheta}$, where for $w = \cdots w_{-1} w_{0} w_{1} \cdots \in X_{\vartheta}$ we let $\vartheta(w)$ denotes the (infinite) set of sequences of the form $v = \cdots v_{-2} v_{-1}.v_0 v_1 \cdots$, with $v_j \in \vartheta(w_j)$ for all $j \in \mathbb{Z}$.
It is easily verified that $\vartheta(X_{\vartheta}) \subset X_{\vartheta}$.

The notation $X_{\vartheta}$ reflects the fact that the random substitution subshift does not depend on the choice of (non-degenerate) probabilities $\bm{P}$.
In fact, this is the case for many structural properties of $\vartheta_{\bm{P}}$.
In these cases, one sometimes refers to $\vartheta$ instead of $\vartheta_{\bm{P}}$ as a random substitution, see for instance \cite{goh2020,grs2019,rus2020,rs2018}.
On the other hand, for some applications, one needs additional structure on the probability space.
In fact, there is an underlying branching process, similar to a Galton--Watson process, that allows one to construct more refined random variables, see \cite{gs2020} for further details.
The measure theoretic properties we consider are typically dependent on the choice of probabilities; however, some of the auxiliary results we use only depend on the set-valued substitution $\vartheta$.
To avoid confusion, for results where there is no dependence on the choice of probabilities we will give the statement in terms of the set-valued substitution $\vartheta$ and omit the dependence on $\bm{P}$ in the notation.

\subsection{Primitive random substitutions}\label{ss:primitivity}
A standard assumption in the study of substitutions (both deterministic and random) is \defn{primitivity}.
Given a random substitution $\vartheta_{\bm{P}} = (\vartheta, \bm{P})$ over an alphabet $\mathcal{A} = \{ a_{1}, \ldots, a_{d} \}$ with cardinality $d \in \mathbb{N}$, we define the \defn{substitution matrix} $M = M_{\vartheta_{\bm{P}}} \in \mathbb{R}^{d \times d}$ of $\vartheta_{\bm{P}}$ by
\begin{align*}
    M_{i, j}
    = \mathbb{E}[\lvert \vartheta_{\bm{P}}  (a_{j}) \rvert_{a_{i}}]
    = \sum_{k = 1}^{r_{j}} p_{j, k} \lvert s^{(j, k)} \rvert_{a_{i}}.
\end{align*}
Since $M$ has only non-negative entries, it has a real eigenvalue $\lambda$ of maximal modulus.
Observe that $\lambda \geq 1$, with $\lambda = 1$ precisely if $M$ is column-stochastic, so that the random substitution is non-expanding.
To avoid this degenerate situation, we assume that $\lambda>1$.
If the matrix $M$ is \defn{primitive} (that is if there exists a $k \in \mathbb{N}$ such that all the entries of $M^{k}$ are positive), the Perron--Frobenius theorem gives that $\lambda$ is a simple eigenvalue and that the corresponding (right) eigenvector $\bm{R} = (R_{1}, \ldots, R_{d})$ can be chosen to have strictly positive entries.
We will normalise this eigenvector so that $\norm{\bm{R}}_1=1$.
We refer to $\lambda$ as the \defn{Perron--Frobenius eigenvalue} of the random substitution, $\vartheta_{\bm{P}}$, with corresponding \defn{Perron--Frobenius eigenvector} $\bm{R}$.
\begin{definition}
    We say that $\vartheta_{\bm{P}}$ is \defn{primitive} if $M = M_{\vartheta_{\bm{P}}}$ is primitive and its Perron--Frobenius eigenvalue satisfies $\lambda > 1$.
\end{definition}

We emphasise that for a random substitution $\vartheta_{\bm{P}}$, being primitive is independent of the (non-degenerate) choice of probabilities $\bm{P}$.
In this sense, primitivity is a property of $\vartheta$ rather than $\vartheta_{\bm{P}}$.

Since $M_{\vartheta_{\bm{P}}}^{k} = M_{\vartheta_{\bm{P}}^{k}}$, the Perron--Frobenius eigenvalue of $\vartheta_{\bm{P}}^{k}$ is $\lambda^{k}$.

\subsection{Compatible random substitutions}\label{ss:compatibility}
Another standard assumption in the study of random substitutions is compatibility, which gives that exact inflation words have a well-defined abelianisation.
In particular, the matrix of a compatible random substitution is independent of the choice of probabilities, so the letter frequencies are uniform and do not depend on the realisation.
As discussed in the introduction, the existence of uniform letter frequencies is fundamental in the proofs of our main results.
\begin{definition}
    We say that a random substitution $\vartheta_{\bm{P}} = (\vartheta, \bm{P})$ is \defn{compatible} if for all $a \in \mathcal{A}$ and $u,v \in \vartheta(a)$, we have $\Phi (u) = \Phi (v)$.
\end{definition}

Compatibility is independent of the choice of probabilities, and a random substitution $\vartheta_{\bm{P}} = (\vartheta, \bm{P})$ is compatible if and only if for all $u \in \mathcal{A}^{+}$, we have that $\lvert s \rvert_{a} = \lvert t \rvert_{a}$ for all $s$ and $t \in \vartheta (u)$, and $a \in \mathcal{A}$.
We write $\lvert \vartheta (u) \rvert_{a}$ to denote this common value, and let $\lvert \vartheta (u) \rvert$ denote the common length of words in $\vartheta (u)$.
For convenience, we write $\lvert \vartheta \rvert = \max_{a \in \mathcal{A}} \lvert \vartheta(a)\rvert$.
For a random substitution that is both primitive and compatible, the (uniform) letter frequencies are encoded by the right Perron--Frobenius eigenvector of the substitution matrix, which by compatibility is independent of the choice of probabilities.
In particular, we have the following (see \cite{que1987} for a proof in the deterministic case, which also holds in the random case by compatibility).
\begin{lemma}[Letter frequency bounds]\label{l:comp-letter-freqs}
    If $\vartheta_{\bm{P}}$ is a primitive and compatible random substitution, then for all $\varepsilon > 0$ there is an integer $N$ such that every word $v$ of length at least $N$ satisfies
    \begin{equation*}
        \lvert v \rvert (R_a - \varepsilon) < \lvert v \rvert_a < \lvert v \rvert (R_a + \varepsilon)
    \end{equation*}
    for all $a \in \mathcal{A}$.
\end{lemma}

The random Fibonacci substitution defined in \cref{ex:random-fib} is compatible, since $\Phi (ab) = \Phi (ba) = (1,1)$.
It is also primitive, since the square of its substitution matrix is positive.
For any choice of probabilities, the right Perron--Frobenius eigenvector is given by $(\tau^{-1},\tau^{-2})$, where $\tau$ denotes the golden ratio.
In terms of letter frequencies, this means that in all sufficiently long legal words, approximately $\tau^{-1}$ proportion of the letters are $a$ and $\tau^{-2}$ proportion are $b$.

The following consequence of \cref{l:comp-letter-freqs} is useful in the proof of \cref{it:lq-inf-equiv}.
\begin{lemma}\label{l:inf-word-sum}
    Let $\vartheta_{\bm{P}} = (\vartheta,\bm{P})$ be a primitive and compatible random substitution and let $q \geq 1$.
    For all $\varepsilon > 0$, there is an $M \in \mathbb{N}$ such that for every $m \geq M$ and $v \in \mathcal{L}_{\vartheta}^{m}$,
    \begin{align*}
        \prod_{a \in \mathcal{A}} \left( \sum_{s \in \vartheta (a)} \mathbb{P} [\vartheta_{\bm{P}} (a) = s]^q \right)^{m (R_a + \varepsilon)}&\leq \sum_{w \in \vartheta (v)} \mathbb{P} [\vartheta_{\bm{P}} (v) = w]^q\\
                                                                                                                                             &\leq \prod_{a \in \mathcal{A}} \left( \sum_{s \in \vartheta (a)} \mathbb{P} [\vartheta_{\bm{P}} (a) = s]^q \right)^{m (R_a - \varepsilon)} \tp
    \end{align*}
    For $q \leq 1$, the same result holds with reversed inequalities.
\end{lemma}
\begin{proof}
    Since $\vartheta_{\bm{P}}$ is compatible, the cutting points of inflation tiles are well-defined, so breaking the sum into inflation tiles we obtain
    \begin{align*}
        \sum_{w \in \vartheta (v)} \mathbb{P} [\vartheta_{\bm{P}} (v) = w]^q &= \sum_{w^1 \in \vartheta (v_1)} \mathbb{P} [\vartheta_{\bm{P}} (v_1) = w]^q \sum_{w^2 \in \vartheta (v_2)} \cdots \sum_{w^m \in \vartheta (v_m)} \mathbb{P} [\vartheta_{\bm{P}} (v_m) = w^m]^q \\
                                                                             &= \prod_{a \in \mathcal{A}} \left( \sum_{s \in \vartheta (a)} \mathbb{P} [\vartheta_{\bm{P}} (a) = s]^q \right)^{\lvert v \rvert_a} \tp
    \end{align*}
    The result then follows by applying \cref{l:comp-letter-freqs} to bound $|v|_a$, noting that for all $a \in \mathcal{A}$ we have $\sum_{s \in \vartheta (a)} \mathbb{P} [\vartheta_{\bm{P}} (a) = s]^q \leq 1$ if $q \geq 1$ and $\sum_{s \in \vartheta (a)} \mathbb{P} [\vartheta_{\bm{P}} (a) = s]^q \geq 1$ if $q \leq 1$.
\end{proof}

\subsection{Frequency measures}\label{ss:freq-measures}
The main object that we associate with a given primitive random substitution $\vartheta_{\bm{P}}$ is the \defn{frequency measure} $\mu_{\bm{P}}$.
This measure quantifies the relative occurrence of a given word in a random substitution.
We now define this measure precisely.

First, we define the \defn{expected frequency} of a word $v \in \mathcal{L}_{\vartheta}$ by
\begin{align*}
    \text{freq} (v)
    = \lim_{k \rightarrow \infty} \frac{\mathbb{E} [ \lvert \vartheta_{\bm{P}}^{k} (a) \rvert_{v} ] }{ \mathbb{E} [\lvert \vartheta_{\bm{P}}^{k} (a)\rvert ]}\tc
\end{align*}
where, by primitivity, this limit is independent of the choice of $a \in \mathcal{A}$.
In fact, we have the stronger property that the word frequencies exist $\mathbb{P}$-almost surely in the limit of large inflation words and are given by $\operatorname{freq}(v)$ for all $v \in \mathcal{L}_{\vartheta}$ (see \cite{gs2020} for further details).
Recalling that $\xi (X_{\vartheta})$ is the algebra of cylinder sets on $X_{\vartheta}$ that specify the origin, we define $\mu_{\bm{P}} \colon \xi (X_{\vartheta}) \cup \{ \varnothing \} \rightarrow [0,1]$ by $\mu_{\bm{P}} (\varnothing) = 0$, $\mu_{\bm{P}}(X_{\vartheta}) = 1$, and $\mu_{\bm{P}} ([v]_{m}) = \operatorname{freq} (v)$ for $v \in \mathcal{L}_{\vartheta}$ and $m \in \{ 1 - \lvert v \rvert, 2 - \lvert v \rvert, \ldots, 0 \}$.
This set function extends to a unique measure (c.f.~\cite[Proposition~5.3 and Theorem~5.9]{gs2020})
\begin{proposition}[\cite{gs2020}]\label{p:freq-measure-exist}
    The set function $\mu_{\bm{P}}$ is a content with mass one which extends uniquely to a shift-invariant ergodic Borel probability measure on $X_\vartheta$.
\end{proposition}

We call the measure $\mu_{\bm{P}}$ defined in \cref{p:freq-measure-exist} the \defn{frequency measure} corresponding to the random substitution $\vartheta_{\bm{P}}$.
Note that frequency measures are dependent on the probabilities of the substitution.
As such, for the subshift of a primitive random substitution that is non-deterministic, there exist uncountably many frequency measures supported on this subshift \cite{gs2020}.
In contrast, the subshift of a primitive deterministic substitution has precisely one frequency measure, which is the unique ergodic measure \cite{que1987}.

Frequency measures corresponding to primitive and compatible random substitutions satisfy the following renormalisation lemma, which relates the measure of a cylinder set of a legal word to measures of cylinder sets of shorter words via the production probabilities of the random substitution.
This result first appeared in \cite{gs2020} and is central to the proof of the main result in \cite{gmrs2023}.
\begin{lemma}[Renormalisation]\label{l:key-lemma}
    Let $\vartheta_{\bm{P}}$ be a primitive and compatible random substitution with corresponding frequency measure $\mu_{\bm{P}}$.
    Let $n \in \mathbb{N}$ and let $k$ be an integer such that every $v \in \mathcal{L}_{\vartheta}^k$ has $|\vartheta (v)| \geq n + |\vartheta (v_1)|$.
    Then for every $u \in \mathcal{L}_{\vartheta}^{n}$,
    \begin{equation*}
        \mu_{\bm{P}}([u]) = \frac{1}{\lambda} \sum_{v \in \mathcal{L}_{\vartheta}^{k}} \mu_{\bm{P}} ([v]) \sum_{j=1}^{\lvert \vartheta (v_1) \rvert} \mathbb{P} [\vartheta_{\bm{P}} (v)_{[j,j+m-1]} = u] \tp
    \end{equation*}
\end{lemma}

\cref{l:key-lemma} is plays an important role in the proof of \cref{it:lq-inf-equiv}, as it relates the sums $\sum_{u \in \mathcal{L}_{\vartheta}^{n}} \mu_{\bm{P}} ([u])$ to sums over smaller words via the production probabilities.
This in turn allows us to obtain relations between $\tau_{\mu_{\bm{P}}}$ and $\varphi_k$.
Under additional assumptions, simplified reformulations of \cref{l:key-lemma} can be obtained (see, for example, \cref{l:recog-key-lemma}, which is used in \cref{it:multi-formalism}).

\subsection{Separation conditions and recognisability}\label{ss:recog}
In this section, we introduce additional common assumptions which either (1) impose a certain separation on inflation words, or (2) impose a certain uniformity on the inflation and the probabilities.
Under these conditions, we can obtain closed-form formulas for the $L^q$-spectrum.
\begin{definition}\label{d:ISC/DSC}
    A random substitution $\vartheta_{\bm{P}} = (\vartheta, \bm{P})$ satisfies the \defn{disjoint set condition} if
    \begin{align*}
        u  \; \text{and} \;  v \in \vartheta (a) \; \text{with} \; u \neq v &\implies \vartheta^{k} (u) \cap \vartheta^{k} (v) = \varnothing
    \end{align*}
    for all $a \in \mathcal{A}$ and $k \in \mathbb{N}$.
    It satisfies the \defn{identical set condition} if
    \begin{align*}
        u \; \text{and} \; v \in \vartheta (a) &\implies \vartheta^{k} (u) = \vartheta^{k} (v)
    \end{align*}
    for all $a \in \mathcal{A}$ and $k \in \mathbb{N}$.
    Moreover, we say that $\vartheta_{\bm{P}}$ has \defn{identical production probabilities} if for all $a \in \mathcal{A}$, $k \in \mathbb{N}$ and $v \in \vartheta^{k} (a)$,
    \begin{align*}
        \mathbb{P} [\vartheta_{\bm{P}}^{k-1} (u_1) = v] = \mathbb{P} [\vartheta_{\bm{P}}^{k-1} (u_2) = v]
    \end{align*}
    for all $u_1$ and $u_2 \in \vartheta (a)$.
\end{definition}

A consequence of the disjoint set condition is that for every $a \in \mathcal{A}$, $k \in \mathbb{N}$ and $w \in \vartheta^k (a)$, there is a unique $v \in \vartheta^{k-1} (a)$ such that $w \in \vartheta (v)$.
In other words, every exact inflation word can be uniquely de-substituted to another exact inflation word.
The following definition extends this idea of unique de-substitution from inflation words to all elements in the subshift.
\begin{definition}\label{d:recog}
    Let $\vartheta_{\bm{P}} = (\vartheta,\bm{P})$ be a primitive and compatible random substitution.
    We call $\vartheta_{\bm{P}}$ \defn{recognisable} if for every $x \in X_{\vartheta}$ there exists a unique $y = \cdots y_{-1}y_{0}y_{1} \cdots \in X_{\vartheta}$ and a unique integer $k \in \{ 0, \ldots, | \vartheta(y_{0}) | - 1 \}$ with $S^{-k}(x) \in \vartheta(y)$.
\end{definition}
The following follows routinely from the definition of recognisability (a proof is given in \cite[Lemma~4.5]{gmrs2023}).
\begin{lemma}\label{l:recog-implies-DSC}
    If $\vartheta_{\bm{P}}$ is a primitive, compatible and recognisable random substitution, then $\vartheta_{\bm{P}}$ satisfies the disjoint set condition.
\end{lemma}
In contrast to the disjoint set condition, recognisability is stable under taking powers (see \cite[Lemma~4.6]{gmrs2023}).
\begin{lemma}\label{l:recognisability-recursion}
    Let $\vartheta_{\bm{P}}$ be a primitive and compatible random substitution and $m \in \mathbb{N}$.
    If $\vartheta_{\bm{P}}$ is recognisable, then so is $\vartheta_{\bm{P}}^m$.
\end{lemma}
An alternative characterisation of recognisability is the following \emph{local} version.
Intuitively, local recognisability means that applying a \emph{finite window} to a sequence is enough to determine the position and the type of the inflation word in the middle of that window.
The following result is given in \cite[Lemma~4.4]{gmrs2023} (see also \cite[Proposition~5.7]{frspreprint}).
\begin{lemma}\label{l:local-recog}
    Let $\vartheta_{\bm{P}} = (\vartheta, \bm{P})$ be a primitive and compatible random substitution.
    If $\vartheta_{\bm{P}}$ is recognisable, then there exists a smallest natural number $\kappa(\vartheta)$, called the \emph{recognisability radius of $\vartheta_{\bm{P}}$}, with the following property: if $x \in \vartheta([a])$ for some $a \in \mathcal{A}$ and $x_{[-\kappa(\vartheta),\kappa(\vartheta)]} = y_{[-\kappa(\vartheta),\kappa(\vartheta)]}$ for some $y \in X_{\vartheta}$, then $y \in \vartheta([a])$.
\end{lemma}
As a consequence of this local characterisation of recognisability, for every legal word $u$ with length greater than twice the radius of recognisability there exists an inflation word $w$, appearing as a subword of $u$, which has a unique decomposition into exact inflation words.
We call the largest such $w$ the \emph{recognisable core} of $u$.

Local recognisability allows us to obtain a stronger version of \cref{l:key-lemma} for recognisable random substitutions.
This result is key to obtaining the coincidence of the $L^q$-spectrum and its inflation word analogue under recognisability for $q<0$, and thus the conclusion of \cref{it:multi-formalism}.
\begin{lemma}\label{l:recog-key-lemma}
    Let $\vartheta_{\bm{P}} = (\vartheta, \bm{P})$ be a primitive and compatible random substitution, with corresponding frequency measure $\mu_{\bm{P}}$ and $u \in \mathcal{L}_{\vartheta}$.
    If $v \in \mathcal{L}_{\vartheta}$ and $w \in \vartheta (v)$ contains $u$ as a subword, then
    \begin{equation*}
        \mu_{\bm{P}} ([u]) \geq \frac{1}{\lambda} \mu_{\bm{P}} ([v]) \mathbb{P} [\vartheta_{\bm{P}} (v) = w] \text{.}
    \end{equation*}
    If, additionally, $\vartheta_{\bm{P}}$ is recognisable, $\lvert u \rvert > 2 \kappa (\vartheta)$ and $w'$ is the recognisable core of $u$ with  $v' \in \mathcal{L}_{\vartheta}$ the unique legal word such that $w' \in \vartheta (v')$, then
    \begin{equation*}
        \mu_{\bm{P}} ([u]) \leq \frac{\kappa(\vartheta)}{\lambda} \mu_{\bm{P}} ([v']) \mathbb{P} [\vartheta_{\bm{P}} (v') = w'] \text{.}
    \end{equation*}
\end{lemma}
\begin{proof}
    If $u$ is a subword of $w \in \vartheta (v)$, then $\mu_{\bm{P}} ([u]) \geq \mu_{\bm{P}} ([w])$.
    Thus applying \cref{l:key-lemma} to $\mu_{\bm{P}} ([w])$,
    \begin{equation*}
        \mu_{\bm{P}} ([u]) \geq \frac{1}{\lambda} \mu_{\bm{P}} ([v]) \mathbb{P} [\vartheta_{\bm{P}} (v) = w] \text{.}
    \end{equation*}
    Now, assume that $\vartheta_{\bm{P}}$ is recognisable, $\lvert u \rvert > 2 \kappa (\vartheta)$ and $w' \in \vartheta (v')$ is the recognisable core of $u$.
    Let $k$ be an integer such that every $t \in \mathcal{L}_{\vartheta}^{k}$ has $\lvert \vartheta (t) \rvert \geq k + \lvert \vartheta (v_1) \rvert$.
    Since there are at most $\kappa (\vartheta)$ letters of $u$ preceding the recognisable core, if $t \in \mathcal{L}_{\vartheta}^{k}$ is a word for which $u \in \vartheta (t)_{[j,j+\lvert u \rvert -1]}$ for some $j \in \{ 1, \ldots , \lvert \vartheta (t_1) \rvert \}$, then $t_i \cdots t_{i + \lvert v \rvert - 1} = v'$ for some $i \in \{ 1, \ldots, \kappa(\vartheta) \}$.
    Moreover, since there is a unique way to decompose $w'$ into exact inflation words, for each $t \in \mathcal{L}_{\vartheta}^{k}$ there can be at most one $j \in \{ 1, \ldots, \vartheta (t_1) \}$ such that $u \in \vartheta (t)_{[j,j+\lvert u \rvert -1]}$.
    Hence, it follows by \cref{l:key-lemma} that
    \begin{align*}
        \mu_{\bm{P}} ([u]) &= \frac{1}{\lambda} \sum_{t \in \mathcal{L}^k} \mu_{\bm{P}} ([t]) \sum_{j=1}^{\lvert \vartheta(t_1) \rvert} \mathbb{P} [\vartheta_{\bm{P}} (t)_{[j,j+\lvert u \rvert-1]} = u] \\
                           &\leq \frac{1}{\lambda} \sum_{i=1}^{\kappa(\vartheta)} \sum_{\substack{t \in \mathcal{L}_{\vartheta}^k \\ t_{i} \cdots t_{i + \lvert v \rvert - 1} = v'}} \mu_{\bm{P}} ([t]) \mathbb{P} [\vartheta_{\bm{P}} (v') = w']\\
                           &= \frac{\kappa(\vartheta)}{\lambda} \mu_{\bm{P}} ([v']) \mathbb{P} [\vartheta_{\bm{P}} (v') = w'] \text{,}
    \end{align*}
    which completes the proof.
\end{proof}

\section{\texorpdfstring{$L^q$}{Lq}-spectra of frequency measures}
In this section, we prove our main results on $L^q$-spectra of frequency measures.
Here, we relate the $L^q$-spectrum to a certain ``symbolic'' $L^q$-spectrum, which we call the \defn{inflation word $L^q$-spectrum}.
Heuristically, the inflation word $L^q$-spectrum is the natural guess for the $L^q$-spectrum if you do not account for non-uniqueness in the positions in which legal words can appear in inflation words.
This notion is introduced in \cref{ss:inf-lq}, where we also state and prove some of its key properties.
In particular, in \cref{p:inflation-isc-dsc}, we prove a simple closed-form formula for the inflation word $L^q$-spectrum under the disjoint set condition or the identical set condition with identical production probabilities.
In \cref{p:mono-bounds}, we establish basic monotonicity results.

Then, in \cref{ss:q>0} and \cref{ss:q<0-lower}, we establish the general bounds for the $L^q$-spectrum in terms of the inflation word $L^q$-spectrum, giving \cref{it:lq-inf-equiv} (the proof is given in \cref{ss:lq-bound-proof}).
We also prove that this bound is sharp in \cref{ss:q<0-recog} under recognisability.
This proves the first part of \cref{it:multi-formalism}.
However this bound need not hold in general: we discuss a counterexample in \cref{ex:q<0-non-inflation}.
Finally, in \cref{ss:en-recov}, we prove differentiability of the $L^q$-spectrum at $q=1$ and show how to recover known results for measure theoretic and topological entropy from our results concerning $L^q$-spectra.

\subsection{Inflation word \texorpdfstring{$L^q$}{Lq}-spectra}\label{ss:inf-lq}
Given a primitive random substitution $\vartheta_{\bm{P}} = (\vartheta, \bm{P})$, we can define an analogue of the $L^q$-spectrum in terms of its production probabilities, in a similar manner to the inflation word analogue of entropy introduced in \cite{gmrs2023}.
In many cases, this notion coincides with the $L^q$-spectrum of the frequency measure associated to $\vartheta_{\bm{P}}$.
For each $k \in \mathbb{N}$ and $q \in \mathbb{R}$, define
\begin{equation*}
    \varphi_k (q) = - \sum_{a \in \mathcal{A}} R_a \log \left( \sum_{s \in \vartheta^k (a)} \mathbb{P} [\vartheta_{\bm{P}}^k (a) = s]^q \right) \tc
\end{equation*}
where $\bm{R}=(R_a)_{a\in\mathcal{A}}$ is the right Perron--Frobenius eigenvector of the substitution matrix of $\vartheta_{\bm{P}}$.
We define the \defn{inflation word $L^q$-spectrum of $\vartheta_{\bm{P}}$} by
\begin{equation*}
    T_{\vartheta,\bm{P}}(q) = \liminf_{k \rightarrow \infty} \frac{\varphi_k(q)}{\lambda^k}\tp
\end{equation*}
We recall that $\lambda$ is the spectral radius of the substitution matrix associated with the random substitution.
We similarly define the upper variant $\lqsymbol{\vartheta,\bm{P}}$ by taking a limit supremum in place of the limit infimum.

We first state some key properties of $T_{\vartheta,\bm{P}} (q)$ which follow easily from the definition.
Firstly, if the random substitution $\vartheta_{\bm{P}}$ is compatible and satisfies either the disjoint set condition or the identical set condition with identical production probabilities, then the limit defining $T_{\vartheta,\bm{P}} (q)$ exists for all $q \in \mathbb{R}$ and is given by a closed-form expression.
For $q \geq 0$, these properties transfer to the $L^q$-spectrum by \cref{it:lq-inf-equiv}.
\begin{proposition}\label{p:inflation-isc-dsc}
    Let $\vartheta_{\bm{P}}$ be a primitive and compatible random substitution and $q\in\R$.
    If $\vartheta_{\bm{P}}$ satisfies the disjoint set condition, then the limit defining $T_{\vartheta,\bm{P}}(q)$ exists and
    \begin{equation*}
        T_{\vartheta,\bm{P}}(q) = \frac{1}{\lambda - 1} \varphi_1 (q)\tp
    \end{equation*}
    If $\vartheta_{\bm{P}}$ satisfies the identical set condition and has identical production probabilities, then the limit defining $T_{\vartheta,\bm{P}}(q)$ exists and
    \begin{equation*}
        T_{\vartheta,\bm{P}}(q) = \frac{1}{\lambda} \varphi_1 (q) \tp
    \end{equation*}
\end{proposition}
\begin{proof}
    Fix $q \in \mathbb{R}$.
    By the Markov property of $\vartheta_{\bm{P}}$, for all $a \in \mathcal{A}$, $k \in \mathbb{N}$ and $v \in \vartheta^k (a)$,
    \begin{equation}\label{eq:prob-relation}
        \mathbb{P} [\vartheta_{\bm{P}}^k (a) = v] = \sum_{s \in \vartheta (a)} \mathbb{P} [\vartheta_{\bm{P}} (a) = s] \, \mathbb{P} [\vartheta_{\bm{P}}^{k-1} (s) = v]\tp
    \end{equation}
    First, suppose $\vartheta_{\bm{P}}$ satisfies the disjoint set condition.
    Then for every $v \in \vartheta^{k} (a)$ there is a unique $s (v) \in \vartheta (a)$ such that $v \in \vartheta^{k-1} (s(v))$.
    Thus, for all $s \in \vartheta(a)$ such that $s \neq s(v)$, we have $\mathbb{P} [\vartheta_{\bm{P}}^{k-1}(s) = v] = 0$, and so it follows by \cref{eq:prob-relation} that
    \begin{align*}
        \sum_{v \in \vartheta^k (a)} \mathbb{P} [\vartheta_{\bm{P}}^k (a) = v]^q &= \sum_{v \in \vartheta^k (a)} \mathbb{P} [\vartheta_{\bm{P}} (a) = s(v)]^q \, \mathbb{P} [\vartheta_{\bm{P}}^{k-1}(s(v)) = v]^q \\
                                                                                     &= \sum_{s \in \vartheta (a)} \mathbb{P} [\vartheta_{\bm{P}} (a) = s]^q \sum_{u \in \vartheta^{k-1} (s)} \mathbb{P} [\vartheta_{\bm{P}}^{k-1} (s) = u]^q \\
                                                                                     &= \left( \sum_{s \in \vartheta (a)} \mathbb{P} [\vartheta_{\bm{P}} (a) = s]^q \right)\prod_{b \in \mathcal{A}} \left( \sum_{u \in \vartheta^{k-1} (b)} \mathbb{P} [\vartheta_{\bm{P}}^{k-1} (b) = u]^q \right)^{\lvert \vartheta (a) \rvert_b}
    \end{align*}
    where in the final equality we use compatibility to split the second sum into inflation tiles.
    Thus
    \begin{align*}
        \varphi_k (q) ={}& - \sum_{a \in \mathcal{A}} R_a \sum_{b \in \mathcal{A}} \lvert \vartheta (a) \rvert_b \log \left( \sum_{u \in \vartheta^{k-1} (b)} \mathbb{P} [\vartheta_{\bm{P}}^{k-1} (b) = u]^q \right)\\
                  &- \sum_{a \in \mathcal{A}} R_a \log \left( \sum_{s \in \vartheta (a)} \mathbb{P} [\vartheta_{\bm{P}} (a) = s]^q \right)\\
        ={}& \lambda \varphi_{k-1} (q) + \varphi_1 (q) \tc
    \end{align*}
    noting that $\sum_{a \in \mathcal{A}} R_a \lvert \vartheta (a) \rvert_b = \lambda R_b$.
    It follows inductively that
    \begin{equation*}
        \frac{1}{\lambda^k} \varphi_k (q) = \sum_{j=1}^{k} \frac{1}{\lambda^j} \varphi_1 (q) \fto{k\to\infty} \frac{1}{\lambda - 1} \varphi_1 (q)\tc
    \end{equation*}
    so the limit defining $T_{\vartheta,\bm{P}} (q)$ exists and is equal to $(\lambda-1)^{-1} \varphi_1 (q)$.

    Next, suppose $\vartheta_{\bm{P}}$ satisfies the identical set condition and has identical production probabilities.
    Then $\mathbb{P} [\vartheta_{\bm{P}}^{k-1} (s^1) = u] = \mathbb{P} [\vartheta_{\bm{P}}^{k-1} (s^2) = u]$ for all $s^1, s^2 \in \vartheta (a)$.
    Hence, it follows by \cref{eq:prob-relation} that
    \begin{equation*}
        \sum_{v \in \vartheta^k (a)} \mathbb{P} [\vartheta_{\bm{P}}^k (a) = v]^q = \sum_{v \in \vartheta^k (a)} \mathbb{P} [\vartheta_{\bm{P}}^{k-1} (s) = v]^q
    \end{equation*}
    for any choice of $s \in \vartheta (a)$.
    By compatibility and the independence of the action,
    \begin{equation*}
        \sum_{v \in \vartheta^k (a)} \mathbb{P} [\vartheta_{\bm{P}}^k (a) = v]^q = \prod_{b \in \mathcal{A}} \left( \sum_{u \in \vartheta^{k-1} (b)} \mathbb{P} [\vartheta_{\bm{P}}^{k-1} (b) = u]^q \right)^{\lvert \vartheta (a) \rvert_b} \tc
    \end{equation*}
    and thus
    \begin{equation*}
        \varphi_k (q) = \sum_{b \in \mathcal{A}} \sum_{a \in \mathcal{A}} R_a \lvert \vartheta (a) \rvert_b \log \left( \sum_{v \in \vartheta^{k-1} (b)} \mathbb{P} [\vartheta_{\bm{P}}^{k-1} (b) = v]^q \right)
        = \lambda \varphi_{k-1} (q) \tc
    \end{equation*}
    noting that $\sum_{a \in \mathcal{A}} R_a \lvert \vartheta (a) \rvert_b = R_b$.
    It follows by induction that $\varphi_k(q)/\lambda^k = \varphi_1(q)/\lambda$ for all $k \in \mathbb{N}$, so we conclude that $T_{\vartheta,\bm{P}}(q)$ exists and equals $\lambda^{-1} \varphi_1 (q)$.
\end{proof}
\begin{proposition}\label{p:mono-bounds}
    Let $\vartheta_{\bm{P}}$ be a primitive and compatible random substitution.
    For all $q>1$ and $q<0$, the sequence $(\lambda^{-k} \varphi_k(q))_k$ is non-decreasing; and for all $0<q<1$, the sequence is non-increasing.
\end{proposition}
\begin{proof}
    This is largely a consequence of Jensen's inequality.
    Note that on the interval $(0,1]$, the function $x \mapsto x^q$ is convex if $q>1$ or $q<0$, and concave if $0<q<1$.
    We first prove this for the case when $q>1$ or $q<0$.
    For all $a \in \mathcal{A}$, $k \in \mathbb{N}$ with $k\geq 2$ and $v \in \vartheta^k (a)$, it follows by the Markov property of $\vartheta_{\bm{P}}$ that
    \begin{align*}
        \sum_{v \in \vartheta^{k}} \mathbb{P} [\vartheta_{\mathbf{P}}^k (a) = v]^q &= \sum_{v \in \vartheta^k (a)} \left( \sum_{s \in \vartheta (a) \colon v \in \vartheta^{k-1} (s)} \mathbb{P} [\vartheta_{\mathbf{P}} (a) = s] \mathbb{P} [\vartheta_{\mathbf{P}}^{k-1} (s) = v] \right)^q\\
                                                                                   &\leq \sum_{v \in \vartheta^k (a)} \left( \frac{ \sum_{s \in \vartheta (a) \colon v \in \vartheta^{k-1} (s)} \mathbb{P} [\vartheta_{\mathbf{P}} (a) = s] \mathbb{P} [\vartheta_{\mathbf{P}}^{k-1} (s) = v]^q }{\sum_{s \in \vartheta (a) \colon v \in \vartheta^{k-1} (s)} \mathbb{P} [\vartheta_{\mathbf{P}} (a) = s]} \right)\\
                                                                                   &\leq\prod_{b \in \mathcal{A}} \left( \sum_{w \in \vartheta^{k-1} (b)} \mathbb{P} [\vartheta_{\mathbf{P}}^{k-1} (b) = w]^q \right)^{\lvert \vartheta (a) \rvert_b} \text{.}
    \end{align*}
    In the second line, we apply Jensen's inequality, and in the third line, we use compatibility to decompose each probability $\mathbb{P} [\vartheta_{\mathbf{P}}^{k-1} (s) = w]$ into inflation tiles.
    It follows that
    \begin{align*}
        \frac{1}{\lambda^k} \varphi_k (q)&\geq - \frac{1}{\lambda^k} \sum_{b \in \mathcal{A}} R_b \sum_{a \in \mathcal{A}} R_a \lvert \vartheta (a) \rvert_b \log \left( \sum_{w \in \vartheta^{k-1} (b)} \mathbb{P} [\vartheta_{\bm{P}}^{k-1} (b) = w]^q \right)\\
                                         &= \frac{1}{\lambda^{k-1}} \varphi_{k-1} (q) \tc
    \end{align*}
    noting that $\sum_{a \in \mathcal{A}} R_a \lvert \vartheta (a) \rvert_b = \lambda$.

    The $0<q<1$ case follows similarly, with Jensen's inequality giving the opposite inequality since $x \mapsto x^q$ is concave.
\end{proof}
An analogous monotonicity result does not hold in general for the $(\lambda^k-1)^{-1} \varphi_k(q)$ bounds, even when $q \geq 0$.
A counterexample is given by the random period doubling substitution (\cref{ex:rpd-end}) with non-uniform probabilities.

\subsection{\texorpdfstring{$L^q$}{Lq}-spectra for non-negative \texorpdfstring{$q$}{q}}\label{ss:q>0}
The majority of the work in proving \cref{it:lq-inf-equiv} lies in proving the bounds in \cref{eq: 0<q<1 bounds}, \cref{eq:q>1 bounds} and \cref{eq:q<0-lower-bounds}.
It suffices to prove the bound for the case $k=1$, since we then obtain the bound for other $k \in \mathbb{N}$ by considering higher powers of the random substitution.
We first prove the upper bound for the case $q>1$.

Throughout this section, we assume that the random substitution is primitive and compatible.
\begin{proposition}\label{p:q>1-upper}
    For all $q > 1$,
    \begin{equation*}
        \overline{\tau}_{\mu_{\bm{P}}} (q) \leq \frac{1}{\lambda - 1} \varphi_{1} (q)\tp
    \end{equation*}
\end{proposition}
\begin{proof}
    Fix $q > 1$.
    Let $\varepsilon > 0$ and, for each $n \in \mathbb{N}$, let $m(n)$ be the integer defined by
    \begin{equation*}
        m(n) = \left\lceil \frac{n}{\lambda - \varepsilon} \right\rceil \tp
    \end{equation*}
    Then the integers $n$ and $m(n)$ satisfy the conditions of \cref{l:key-lemma}, so it follows that
    \begin{equation*}
        \sum_{u \in \mathcal{L}_{\vartheta}^{n}} \mu_{\bm{P}}([u])^q = \sum_{u \in \mathcal{L}_{\vartheta}^{n}} \left( \frac{1}{\lambda} \sum_{v \in \mathcal{L}_{\vartheta}^{m(n)}} \mu_{\bm{P}}([v]) \sum_{j=1}^{\lvert \vartheta (v_1) \rvert} \mathbb{P} [\vartheta_{\bm{P}} (v)_{[j,j+n-1]} = u] \right)^q \tp
    \end{equation*}
    Since $q > 1$, the function $x \mapsto x^q$ is superadditive on the interval $[0,1]$, so
    \begin{align*}
        \sum_{u \in \mathcal{L}_{\vartheta}^{n}} \mu_{\bm{P}}([u])^q &\geq \sum_{u \in \mathcal{L}_{\vartheta}^{n}} \sum_{v \in \mathcal{L}_{\vartheta}^{m(n)}} \mu_{\bm{P}}([v])^q \left( \frac{1}{\lambda} \sum_{j=1}^{\lvert \vartheta (v_1) \rvert} \mathbb{P} [\vartheta_{\bm{P}} (v)_{[j,j+n-1]} = u] \right)^q\\
        &\geq \frac{1}{\lambda^q} \sum_{v \in \mathcal{L}_{\vartheta}^{m(n)}} \mu_{\bm{P}}([v])^q \sum_{j=1}^{\lvert \vartheta (v_1) \rvert} \sum_{u \in \mathcal{L}_{\vartheta}^{n}} \mathbb{P} [\vartheta_{\bm{P}} (v)_{[j,j+n-1]} = u]^q \tp
    \end{align*}
    We now bound the probability on the right of this expression by the production probability of an inflation word.
    If $w(u) \in \vartheta (v)$ contains $u$ as a subword in position $j$, then $\mathbb{P} [\vartheta_{\bm{P}} (v)_{[j,j+n-1]} = u] \geq \mathbb{P} [\vartheta_{\bm{P}} (v) = w(u)]$.
    Hence,
    \begin{equation*}
        \sum_{u \in \mathcal{L}_{\vartheta}^{n}} \mathbb{P} [\vartheta_{\bm{P}} (v)_{[j, j + n - 1]} = u]^q \geq \sum_{w \in \vartheta (v)} \mathbb{P} [\vartheta_{\bm{P}} (v) = w]^q
    \end{equation*}
    for all $j \in \{1, \ldots, \lvert \vartheta (v_1) \rvert \}$.

    Since $\vartheta_{\bm{P}}$ is compatible, by \cref{l:inf-word-sum} there exists an $N \in \mathbb{N}$ such that for all $n \geq N$ and all $v \in \mathcal{L}_{\vartheta}^{m(n)}$
    \begin{equation*}
        \sum_{w\in\vartheta(v)}\mathbb{P} [\vartheta_{\bm{P}} (v) = w]^q \geq \prod_{a \in \mathcal{A}} \left( \sum_{s \in \vartheta (a)} \mathbb{P} [\vartheta_{\bm{P}} (a) = s]^q \right)^{m(n) (R_a + \varepsilon)} \tp
    \end{equation*}
    Hence,
    \begin{equation*}
        \sum_{u \in \mathcal{L}_{\vartheta}^{n}} \mu_{\bm{P}}([u])^q \geq \frac{1}{\lambda^q} \prod_{a \in \mathcal{A}} \left( \sum_{s \in \vartheta (a)} \mathbb{P} [\vartheta_{\bm{P}} (a) = s]^q \right)^{m(n) (R_a + \varepsilon)} \sum_{v \in \mathcal{L}_{\vartheta}^{m(n)}} \mu_{\bm{P}}([v])^q \tp
    \end{equation*}
    Taking logarithms, rearranging and dividing by $n$ gives
    \begin{align*}
        - \frac{1}{n} \log \left( \sum_{u \in \mathcal{L}_{\vartheta}^{n}} \mu_{\bm{P}}([u])^q \right) \leq{}& - \frac{1}{n} \log \left( \sum_{v \in \mathcal{L}_{\vartheta}^{m(n)}} \mu_{\bm{P}}([v])^q \right) + \frac{1}{n} \log \lambda^q\\
        &- \frac{m(n)}{n} \sum_{a \in \mathcal{A}} (R_a + \varepsilon) \log \left( \sum_{s \in \vartheta (a)} \mathbb{P} [\vartheta_{\bm{P}} (a) = s]^q \right) \tp
    \end{align*}
    Noting that $m(n)/n \rightarrow (\lambda-\varepsilon)^{-1}$ as $n \rightarrow \infty$, it follows by \cref{l:lq-subseq} that
    \begin{equation*}
        \overline{\tau}_{\mu_{\bm{P}}} (q) \leq \frac{1}{\lambda - \varepsilon} \overline{\tau}_{\mu_{\bm{P}}} (q) + \frac{1}{\lambda - \varepsilon} \sum_{a \in \mathcal{A}} \log \left( \sum_{s \in \vartheta (a)} \mathbb{P} [\vartheta_{\bm{P}} (a) = s]^q \right) + c \varepsilon
    \end{equation*}
    where $c\coloneqq (\# \mathcal{A}) \max_{a \in \mathcal{A}} \log (\sum_{s \in \vartheta (a)} \mathbb{P} [\vartheta_{\bm{P}} (a) = s]^q)$.
    But $\varepsilon>0$ was arbitrary; letting $\varepsilon \rightarrow 0$ and rearranging we obtain
    \begin{equation*}
        \overline{\tau}_{\mu_{\bm{P}}} (q) \leq \frac{1}{\lambda - 1} \varphi_1 (q) \tc
    \end{equation*}
    which completes the proof.
\end{proof}

We now prove the corresponding lower bound.
\begin{proposition}\label{p:q>1-lower}
    For all $q > 1$,
    \begin{equation*}
        \tau_{\mu_{\bm{P}}} (q) \geq \frac{1}{\lambda} \varphi_{1} (q) \tp
    \end{equation*}
\end{proposition}
\begin{proof}
    Let $\varepsilon > 0$ and, for each $n \in \mathbb{N}$, let $m(n)$ be the integer defined by
    \begin{equation*}
        m(n) = \left\lceil \frac{n}{\lambda - \varepsilon} \right\rceil \tp
    \end{equation*}
    Since $q>1$, the function $x \mapsto x^q$ is convex on the interval $[0,1]$.
    Hence, it follows by \cref{l:key-lemma} and two applications of Jensen's inequality that
    \begin{align*}
        \sum_{u \in \mathcal{L}_{\vartheta}^{n}} \mu_{\bm{P}}([u])^q &= \sum_{u \in \mathcal{L}_{\vartheta}^{n}} \left( \frac{1}{\lambda} \sum_{v \in \mathcal{L}_{\vartheta}^{m(n)}} \mu_{\bm{P}}([v]) \sum_{j=1}^{\lvert \vartheta (v_1) \rvert} \mathbb{P} [\vartheta_{\bm{P}} (v)_{[j,j+n-1]} = u] \right)^q\\
        &\leq \sum_{v \in \mathcal{L}_{\vartheta}^{m(n)}} \mu_{\bm{P}}([v]) \sum_{u \in \mathcal{L}_{\vartheta}^{n}} \left( \frac{1}{\lambda} \sum_{j=1}^{\lvert \vartheta (v_1) \rvert} \mathbb{P} [\vartheta_{\bm{P}} (v)_{[j,j+n-1]} = u] \right)^q\\
        &\leq \frac{|\vartheta|^{q-1}}{\lambda^q} \sum_{v \in \mathcal{L}_{\vartheta}^{m(n)}} \mu_{\bm{P}}([v]) \sum_{j=1}^{\lvert \vartheta (v_1) \rvert} \sum_{u \in \mathcal{L}_{\vartheta}^{n}} \mathbb{P} [\vartheta_{\bm{P}} (v)_{[j,j+n-1]} = u]^q\tp
    \end{align*}
    We bound above the probability on the right of this expression by the production probability of a sufficiently large inflation word contained in $u$.
    By compatibility, there is an integer $k(n)$ such that $j+n \leq \lvert \vartheta (v_{[1,m(n)-k(n)]})\rvert$ for all $n \in \mathbb{N}$ and $v \in \mathcal{L}_{\vartheta}^{m(n)}$, where $\lim k(n)/n=0$.
    In particular, for every $v \in \mathcal{L}_{\vartheta}^{n}$, a realisation of $\vartheta (v_{[2,m(n)-k(n)]})$ is contained in $u$ as an inflation word, so
    \begin{equation*}
        \sum_{u \in \mathcal{L}_{\vartheta}^{n}} \mathbb{P} [\vartheta_{\bm{P}} (v)_{[j,j+ n - 1]} = u]^q
        \leq \sum_{w \in \vartheta (v_2 \cdots v_{m(n)-k(n)})} \mathbb{P} [\vartheta_{\bm{P}} (v_2 \cdots v_{m(n)-k(n)}) = w]^q \tp
    \end{equation*}
    We now bound this quantity uniformly for all $v \in \mathcal{L}_{\vartheta}^{m(n)}$.
    By \cref{l:inf-word-sum} and the above, there is an $N \in \mathbb{N}$ such that for all $n \geq N$
    \begin{equation*}
        \sum_{u \in \mathcal{L}_{\vartheta}^{n}} \mu_{\bm{P}}([u])^q \leq \frac{|\vartheta|^{q-1}}{\lambda^q}\prod_{a \in \mathcal{A}} \left( \sum_{s \in \vartheta (a)} \mathbb{P} [\vartheta_{\bm{P}} (a) = s]^q \right)^{(m(n)-k(n)-1) (R_a - \varepsilon)} \tp
    \end{equation*}
    Taking logarithms, rearranging and dividing by $n$ gives
    \begin{align*}
        - \frac{1}{n} \log \left( \sum_{u \in \mathcal{L}_{\vartheta}^{n}} \mu_{\bm{P}}([u])^q\!\right) &\geq \frac{m(n)-k(n)-1}{n} \sum_{a \in \mathcal{A}} (R_a - \varepsilon) \log \left( \sum_{s \in \vartheta (a)} \mathbb{P} [\vartheta_{\bm{P}} (a) = s]^q \right)\\
                                                                                                       &\qquad-\frac{\log(|\vartheta|^{q-1}/\lambda^q)}{n}\\
                                                                                                       &\xrightarrow{n \rightarrow \infty} \frac{1}{\lambda- \varepsilon} \sum_{a \in \mathcal{A}} (R_a - \varepsilon) \log \left( \sum_{s \in \vartheta (a)} \mathbb{P} [\vartheta_{\bm{P}} (a) = s]^q \right) \tc
    \end{align*}
    But $\varepsilon>0$ was arbitrarily, so
    \begin{equation*}
        \tau_{\mu_{\bm{P}}} (q) \geq \frac{1}{\lambda} \varphi_1 (q) \tc
    \end{equation*}
    which completes the proof.
\end{proof}

We now state the bounds for the $q \in (0,1)$ case.
We do not give a proof here since the arguments mirror the proofs of \cref{p:q>1-upper} and \cref{p:q>1-lower}, except with reversed inequalities since $x\mapsto x^q$ is concave rather than convex and subadditive as opposed to superadditive.
\begin{proposition}\label{p:0<q<1-bounds}
    If $q \in (0,1)$, then
    \begin{equation*}
        \frac{1}{\lambda-1} \varphi_1 (q) \leq \tau_{\mu_{\bm{P}}} (q) \leq \lqupper{\mu_{\bm{P}}}(q) \leq \frac{1}{\lambda} \varphi_1 (q) \tp
    \end{equation*}
\end{proposition}

\subsection{\texorpdfstring{$L^q$}{Lq}-spectra for negative \texorpdfstring{$q$}{q}: lower bounds}\label{ss:q<0-lower}
For $q<0$, there exist primitive and compatible random substitutions for which $\tau_{\mu_{\bm{P}}} (q)$ and $T_{\vartheta,\bm{P}} (q)$ do not coincide (see, for instance, \cref{ex:q<0-non-inflation}).
However, we still obtain that $\tau_{\mu_{\bm{P}}} (q) \geq T_{\vartheta,\bm{P}} (q)$ for all $q<0$.
To prove this, it suffices to show the sequence of bounds in \cref{eq:q<0-lower-bounds} holds.
Again, we only need to prove the bound for $k=1$ since the remaining bounds follow by considering powers of the random substitution.
\begin{proposition}\label{p:q<0-lower}
    If $\vartheta_{\bm{P}}$ is a primitive and compatible random substitution, then for all $q < 0$,
    \begin{equation*}
        \tau_{\mu_{\bm{P}}} (q) \geq \frac{1}{\lambda - 1} \varphi_{1} (q) \tp
    \end{equation*}
\end{proposition}
\begin{proof}
    Let $\varepsilon > 0$ be sufficiently small and for $n$ sufficiently large, let $m(n)$ be the integer defined by
    \begin{equation*}
        m(n) = \left\lceil \frac{n}{\lambda - \varepsilon} \right\rceil \tp
    \end{equation*}
    To avoid division by zero, we rewrite \cref{l:key-lemma} in a form where we do not sum over elements equal to zero.
    Here, we write $u \blacktriangleleft \vartheta (v)$ to mean there is a realisation $w$ of $\vartheta (v)$ for which $u$ appears as a subword of $w$.
    For each $v \in \mathcal{L}_{\vartheta}^{m(n)}$ and $u \in \mathcal{L}_{\vartheta}^{n}$, let $\mathcal{J} (v,u) = \{ j \in \{ 1,\ldots, \lvert \vartheta (v_1) \rvert \} : u \in \vartheta (v)_{[j,j+n-1]} \}$.
    If $j \notin \mathcal{J} (u,v)$, then $\mathbb{P} [\vartheta_{\bm{P}} (v)_{[j,j+n-1]} = u] = 0$, and if $u$ does not appear as a subword of any realisations of $\vartheta (v)$, then $\mathcal{J} (u,v) = \varnothing$.
    Therefore, we can rewrite \cref{l:key-lemma} as
    \begin{equation*}
        \mu_{\bm{P}}([u]) = \frac{1}{\lambda} \sum_{\substack{v \in \mathcal{L}_{\vartheta}^{m(n)} \\ u \blacktriangleleft \vartheta (v)}} \mu_{\bm{P}}([v]) \sum_{j \in \mathcal{J} (v,u)} \mathbb{P} [\vartheta_{\bm{P}} (v)_{[j,j+n-1]} = u] \tp
    \end{equation*}
    Hence, by the subadditivity of the function $x \mapsto x^q$ on the domain $(0,1]$,
    \begin{align*}
        \sum_{u \in \mathcal{L}_{\vartheta}^{n}} \mu_{\bm{P}}([u])^q &= \sum_{u \in \mathcal{L}_{\vartheta}^{n}} \left( \frac{1}{\lambda} \sum_{\substack{v \in \mathcal{L}_{\vartheta}^{m(n)} \\ u \blacktriangleleft \vartheta (v)}} \mu_{\bm{P}}([v]) \sum_{j \in \mathcal{J} (v,u)} \mathbb{P} [\vartheta_{\bm{P}} (v)_{[j,j+n-1]} = u] \right)^q \\
        &\leq \frac{1}{\lambda^q} \sum_{u \in \mathcal{L}_{\vartheta}^{n}} \sum_{\substack{v \in \mathcal{L}_{\vartheta}^{m(n)} \\ u \blacktriangleleft \vartheta (v)}} \mu_{\bm{P}}([v])^q \sum_{j \in \mathcal{J} (v,u)} \mathbb{P} [\vartheta_{\bm{P}} (v)_{[j,j+n-1]} = u]^q\\
        &= \frac{1}{\lambda^q} \sum_{v \in \mathcal{L}_{\vartheta}^{m(n)}} \mu_{\bm{P}}([v])^q \sum_{\substack{u \in \mathcal{L}_{\vartheta}^{n} \\ u \blacktriangleleft \vartheta (v)}} \sum_{j \in \mathcal{J} (v,u)} \mathbb{P} [\vartheta_{\bm{P}} (v)_{[j,j+n-1]} = u]^q \tp
    \end{align*}
    For each $j \in \mathcal{J} (v,u)$, let $w_{j} (u) \in \vartheta (v)$ be a word such that $w_j (u)_{[j,j+n-1]} = u$.
    Note that there are at most $K\coloneqq 2|\vartheta|(\# \mathcal{A})^{\lvert \vartheta \rvert}$ different $u \in \mathcal{L}_{\vartheta}^{n}$ such that $w_j (u)_{[j,j+n-1]} = u$.
    Hence,
    \begin{align*}
        \sum_{\substack{u \in \mathcal{L}_{\vartheta}^{n} \\ u \blacktriangleleft \vartheta (v)}} \sum_{j \in \mathcal{J} (v,u)} \mathbb{P} [\vartheta_{\bm{P}} (v)_{[j,j+n-1]} = u]^q&\leq \sum_{\substack{u \in \mathcal{L}_{\vartheta}^{n} \\ u \blacktriangleleft \vartheta (v)}} \sum_{j \in \mathcal{J} (v,u)} \mathbb{P} [\vartheta_{\bm{P}} (v) = w_j (u)]^q\\
        &\leq K\sum_{w \in \vartheta (v)} \mathbb{P} [\vartheta_{\bm{P}} (v) = w]^q
    \end{align*}
    and it follows that
    \begin{equation*}
        \sum_{u \in \mathcal{L}_{\vartheta}^{n}} \mu_{\bm{P}}([u])^q \leq \lambda^{-q}K\sum_{v \in \mathcal{L}_{\vartheta}^{m(n)}} \mu_{\bm{P}}([v])^q \sum_{w \in \vartheta (v)} \mathbb{P} [\vartheta_{\bm{P}} (v) = w]^q \tp
    \end{equation*}
    Thus, by \cref{l:inf-word-sum}, for all $\varepsilon > 0$ there is an integer $N$ such that for all $n \geq N$
    \begin{equation*}
        \sum_{u \in \mathcal{L}_{\vartheta}^{n}} \mu_{\bm{P}}([u])^q \leq \lambda^{-q}K\prod_{a \in \mathcal{A}} \left( \sum_{s \in \vartheta (a)} \mathbb{P} [\vartheta_{\bm{P}} (a) = s]^q \right)^{m(n)(R_a + \varepsilon)}  \left( \sum_{v \in \mathcal{L}_{\vartheta}^{m(n)}} \mu_{\bm{P}}([v])^q \right) \tp
    \end{equation*}
    Taking logarithms, rearranging and dividing by $n$ gives
    \begin{align*}
        - \frac{1}{n} \log \left( \sum_{u \in \mathcal{L}_{\vartheta}^{n}} \mu_{\bm{P}}([u])^q \right) \geq{}& - \frac{1}{n} \log \left( \sum_{v \in \mathcal{L}_{\vartheta}^{m(n)}} \mu_{\bm{P}}([v])^q \right) + \frac{1}{n} \log ( \lambda^{-q}K)\\
        &- \frac{m(n)}{n} \sum_{a \in \mathcal{A}} (R_a + \varepsilon) \log \left( \sum_{s \in \vartheta (a)} \mathbb{P} [\vartheta_{\bm{P}} (a) = s]^q \right) \tp
    \end{align*}
    Noting that $m(n)/n \rightarrow (\lambda-\varepsilon)^{-1}$ as $n \rightarrow \infty$, it follows by \cref{l:lq-subseq} that
    \begin{equation*}
        \tau_{\mu_{\bm{P}}} (q) \geq \frac{1}{\lambda - \varepsilon} \tau_{\mu_{\bm{P}}} (q) + \frac{1}{\lambda - \varepsilon} \sum_{a \in \mathcal{A}} \log \left( \sum_{s \in \vartheta (a)} \mathbb{P} [\vartheta_{\bm{P}} (a) = s]^q \right) + c \varepsilon
    \end{equation*}
    where $c\coloneqq (\#\mathcal{A}) \max_{a \in \mathcal{A}} \log (\sum_{s \in \vartheta (a)} \mathbb{P} [\vartheta_{\bm{P}} (a) = s]^q)$.
    Letting $\varepsilon \rightarrow 0$ and rearranging, we obtain
    \begin{equation*}
        \overline{\tau}_{\mu_{\bm{P}}} (q) \geq \frac{1}{\lambda - 1} \varphi_1 (q) \tc
    \end{equation*}
    which completes the proof.
\end{proof}

\subsection{Proof of general bounds for the \texorpdfstring{$L^q$}{Lq} spectrum}\label{ss:lq-bound-proof}
Using the bounds proved in the prior two sections, we can now complete the proof of \cref{it:lq-inf-equiv}.
\begin{proof}[of \cref{it:lq-inf-equiv}]
    Since, for each $k \in \mathbb{N}$, the random substitution $\vartheta_{\bm{P}}^{k}$ gives rise to the same frequency measure as $\vartheta_{\bm{P}}$, applying \cref{p:q>1-upper}, \cref{p:q>1-lower} and \cref{p:0<q<1-bounds} to $\vartheta_{\bm{P}}^k$,
    \begin{equation*}
        \frac{1}{\lambda^k} \varphi_k (q) \leq \tau_{\mu_{\bm{P}}} (q) \leq \overline{\tau}_{\mu_{\bm{P}}} (q) \leq \frac{1}{\lambda^k - 1} \varphi_k (q)
    \end{equation*}
    for all $q > 1$ and
    \begin{equation*}
        \frac{1}{\lambda^k - 1} \varphi_k (q) \leq \tau_{\mu_{\bm{P}}} (q) \leq \overline{\tau}_{\mu_{\bm{P}}} (q) \leq \frac{1}{\lambda^k} \varphi_k (q)
    \end{equation*}
    for $0<q<1$.
    Letting $k \rightarrow \infty$ gives
    \begin{equation*}
        \tau_{\mu_{\bm{P}}} (q) = \lqupper{\mu_{\bm{P}}} (q) = T_{\vartheta,\bm{P}} (q) = \lqsymbol{\vartheta,\bm{P}}(q)
    \end{equation*}
    for all $q \in (0,1) \cup (1, \infty)$, so the limits defining $\tau_{\mu_{\bm{P}}} (q)$ and $T_{\vartheta,\bm{P}}(q)$ both exist and coincide.
    The same holds for $q=0$ and $q=1$ by continuity.
    The monotonicity of the bounds $\lambda^{-k} \varphi_k (q)$ follows by \cref{p:mono-bounds}.
    Finally for $q<0$, for each $k \in \mathbb{N}$, applying \cref{p:q<0-lower} to $\vartheta_{\bm{P}}^k$ gives that
    \begin{equation*}
        \tau_{\mu_{\bm{P}}} (q) \geq \frac{1}{\lambda^k-1} \varphi_k (q) \tp
    \end{equation*}
    Passing to the limit completes the proof.
\end{proof}

\subsection{\texorpdfstring{$L^q$}{Lq}-spectra for negative \texorpdfstring{$q$}{q} under recognisability}\label{ss:q<0-recog}
While the upper bound does not hold in general for $q < 0$, for recognisable random substitutions we can obtain this using \cref{l:recog-key-lemma}, which we recall is a refinement of \cref{l:key-lemma} using recognisability.
\begin{proposition}\label{p:q<0-lower-recognisable}
    If $\vartheta_{\bm{P}}$ is a primitive, compatible and recognisable random substitution, then for all $q<0$,
    \begin{equation*}
        \overline{\tau}_{\mu_{\bm{P}}} (q) \leq \frac{1}{\lambda - 1} \varphi_1 (q) \tp
    \end{equation*}
\end{proposition}
\begin{proof}
    Let $\varepsilon > 0$ be sufficiently small and, for each $n \in \mathbb{N}$ sufficiently large, let $m(n)$ be the integer defined by
    \begin{equation*}
        m(n) = \left\lfloor \frac{n}{\lambda - \varepsilon} \right\rfloor \tp
    \end{equation*}
    For each $u \in \mathcal{L}_{\vartheta}^{n + 2 \kappa (\vartheta)}$, let $w(u)$ denote the recognisable core of $u$.
    Further, let $v(u)$ denote the unique legal word such that $w(u) \in \vartheta (v(u))$.
    Then, by \cref{l:recog-key-lemma},
    \begin{equation}\label{eq: q<0 recog proof eq1}
        \mu_{\bm{P}} ([u]) \leq \frac{\kappa (\vartheta)}{\lambda} \mu_{\bm{P}} ([v(u)]) \mathbb{P} [\vartheta_{\bm{P}}(v(u)) = w(u)] \text{.}
    \end{equation}
    For all $u \in \mathcal{L}_{\vartheta}^{n + 2 \kappa (\vartheta)}$, the recognisable core $w(u)$ has length at least $n$ so, by compatibility, there is an integer $N$ such that if $n \geq N$, then $\lvert v(u) \rvert \geq m(n)$ for all $u \in \mathcal{L}_{\vartheta}^{n + 2 \kappa (\vartheta)}$.
    In particular, for every $u$ there exists a $v \in \mathcal{L}_{\vartheta}^{m(n)}$ such that $\mu_{\bm{P}} ([v(u)]) \leq \mu_{\bm{P}} ([v])$ and a $w \in \vartheta (v)$ such that $\mathbb{P} [\vartheta_{\bm{P}} (v(u)) = w(u)] \leq \mathbb{P} [\vartheta_{\bm{P}} (v) = w]$.
    Hence, it follows by \eqref{eq: q<0 recog proof eq1} and \cref{l:inf-word-sum} that
    \begin{align*}
        \sum_{u \in \mathcal{L}_{\vartheta}^{n + 2 \kappa (\vartheta)}} \mu_{\bm{P}}([u])^q &\geq \frac{1}{\lambda^q} \sum_{v \in \mathcal{L}_{\vartheta}^{m(n)}} \mu_{\bm{P}}([v])^q \sum_{w \in \vartheta (v)} \mathbb{P} [\vartheta_{\bm{P}} (v) = w]^q\\
                                                                                            &\geq \frac{1}{\lambda^q} \prod_{a \in \mathcal{A}} \left( \sum_{s \in \vartheta (a)} \mathbb{P} [\vartheta_{\bm{P}} (a) = s]^q \right)^{m (R_a - \varepsilon)} \sum_{v \in \mathcal{L}_{\vartheta}^{m(n)}} \mu_{\bm{P}}([v])^q \tc
    \end{align*}
    noting that since $q < 0$, the function $x \mapsto x^q$ is decreasing on $(0,1]$.
    Taking logarithms, rearranging and dividing by $n$ gives
    \begin{align*}
        - \frac{1}{n} \log \left( \sum_{u \in \mathcal{L}_{\vartheta}^{n}} \mu_{\bm{P}}([u])^q \right) \leq{}& - \frac{1}{n} \log \left( \sum_{v \in \mathcal{L}_{\vartheta}^{m(n)}} \mu_{\bm{P}}([v])^q \right) + \frac{1}{n} \log \lambda^q\\
                                                                                                             &- \frac{m(n)}{n} \sum_{a \in \mathcal{A}} (R_a - \varepsilon) \log \left( \sum_{s \in \vartheta (a)} \mathbb{P} [\vartheta_{\bm{P}} (a) = s]^q \right) \tp
    \end{align*}
    Noting that $m(n)/n \rightarrow (\lambda - \varepsilon)^{-1}$ as $n \rightarrow \infty$, it follows by \cref{l:lq-subseq} that
    \begin{equation*}
        \overline{\tau}_{\mu_{\bm{P}}} (q) \leq \frac{1}{\lambda - \varepsilon} \overline{\tau}_{\mu_{\bm{P}}} (q) + \frac{1}{\lambda - \varepsilon} \sum_{a \in \mathcal{A}} \log \left( \sum_{s \in \vartheta (a)} \mathbb{P} [\vartheta_{\bm{P}} (a) = s]^q \right) + c \varepsilon
    \end{equation*}
    where $c\coloneqq (\# \mathcal{A}) \max_{a \in \mathcal{A}} \log (\sum_{s \in \vartheta (a)} \mathbb{P} [\vartheta_{\bm{P}} (a) = s]^q)$.
    Letting $\varepsilon \rightarrow 0$ and rearranging, we obtain
    \begin{equation*}
        \overline{\tau}_{\mu_{\bm{P}}} (q) \leq \frac{1}{\lambda - 1} \varphi_1 (q) \tc
    \end{equation*}
    which completes the proof.
\end{proof}

\subsection{Recovering entropy from the \texorpdfstring{$L^q$}{Lq}-spectrum}\label{ss:en-recov}
Since the $L^q$-spectrum encodes both topological and measure theoretic entropy, \cref{it:lq-inf-equiv} provides an alternative means of proving the coincidence of these quantities with the inflation word analogues introduced in \cite{goh2020, gmrs2023} for primitive and compatible random substitutions.

For notational simplicity, set
\begin{equation*}
    \rho_k=-\sum_{a \in \mathcal{A}} R_a \sum_{s \in \vartheta^k (a)}\mathbb{P}[\vartheta_{\bm{P}}^k (a) = s]\log(\mathbb{P}[\vartheta_{\bm{P}}^k (a) = s])\tp
\end{equation*}
\begin{proof}[of \cref{ic:entropy-recover}]
    We first establish the result for topological entropy.
    By \cref{it:lq-inf-equiv}, the limit defining $T_{\vartheta,\bm{P}} (0)$ exists; in particular,
    \begin{equation*}
        \lim_{m \rightarrow \infty} \frac{1}{\lambda^k} \sum_{a \in \mathcal{A}} R_a \log (\# \vartheta^m (a))
    \end{equation*}
    exists.
    Since $\htop(X_{\vartheta}) = - \tau_{\mu_{\bm{P}}} (0) = - T_{\vartheta,\bm{P}} (0)$, we conclude that
    \begin{equation*}
        \htop(X_{\vartheta}) = - \lim_{m \rightarrow \infty} \frac{1}{\lambda^k} \sum_{a \in \mathcal{A}} R_a \log (\# \vartheta^m (a))
    \end{equation*}
    as claimed.

    Now we consider measure theoretic entropy.
    We first make the following elementary observation: if $f$ and $g$ are concave functions with $f(1)=g(1)$ and $f(x)\leq g(x)$ for all $x\geq 1$, then $f^+(1)\leq g^+(1)$.
    Indeed, for all $\epsilon>0$,
    \begin{equation*}
        \frac{f(1+\epsilon)-f(1)}{\epsilon}\leq\frac{g(1+\epsilon)-g(1)}{\epsilon}\tc
    \end{equation*}
    and taking the limit as $\epsilon$ goes to $0$ (which always exists by concavity) yields the desired inequality.

    Recall that $\tau_{\mu_{\bm{P}}}$ and $\lambda^{-k}\varphi_k$ are concave functions with $\tau_{\mu_{\bm{P}}}(1)=\varphi_k(1)=0$ for all $k\in\N$.
    Moreover, $\varphi_k$ is differentiable for all $k\in\N$ with $\varphi_k'(1)=\rho_k$ and by \cref{p:mono-bounds} and \cref{it:lq-inf-equiv}, $\bigl(\lambda^{-k}\varphi_k\bigr)_{j=1}^\infty$ converges monotonically to $\tau_{\mu_{\bm{P}}}$ from below.
    In particular, $\rho_k/\lambda^k$ is a monotonically increasing sequence bounded above by $\tau_{\mu_{\bm{P}}}^+(1)$, so that the limit indeed exists.
    Thus
    \begin{equation*}
        \tau_{\mu_{\bm{P}}}^+(1)=\lim_{k\to\infty}\frac{\rho_k}{\lambda^k}
    \end{equation*}
    since $\varphi_k(q)/(\lambda^k-1)\geq\tau_{\mu_{\bm{P}}}(q)$ for all $q\in(0,\infty)$, using the preceding observation.

    The result for $\tau_{\mu_{\bm{P}}}^-(1)$ follows by an identical argument, instead using monotonicity and the corresponding bounds for $q\in(0,1)$.
    Thus $\tau_{\mu_{\bm{P}}}'(1)=\lim_{k\to\infty}\rho_k/\lambda^k$, so the desired result follows by \cref{l:lq-core}(c).
\end{proof}

\section{Recognisability and the multifractal formalism}
In this section we establish the multifractal formalism as stated in \cref{it:multi-formalism}.
Our proof will follow from a \emph{constrained variational principle}, which is obtained by considering typical local dimensions of one frequency measure $\mu_{\bm{P}}$ relative to another frequency measure $\mu_{\bm{Q}}$.
Our strategy is to prove the almost sure existence of \emph{relative letter frequencies} in \cref{l:expected-inf-frequencies}: this result, combined with recognisability, gives \cref{p:relative-local-dim}.
The multifractal formalism then follows from this dimensional result combined with the formula for the $L^q$-spectrum proved in \cref{p:q<0-lower-recognisable}---the proof is given in \cref{ss:multi-proof}.

\subsection{Non-typical local dimensions}
To prove the multifractal formalism for a given frequency measure $\mu_{\bm{P}}$, we show that for every $\alpha \in[\alpha_{\min}, \alpha_{\max}]$, there exists another frequency measure $\mu_{\bm{Q}}$ such that $\dimH \mu_{\bm{Q}} \geq \tau_{\mu_{\bm{P}}}^{*} (\alpha)$ and $\ldim (\mu_{\bm{P}}, x) = \alpha$ for $\mu_{\bm{Q}}$-almost every $x \in X_{\vartheta}$.
Given a primitive set-valued substitution $\vartheta$, permissible probabilities $\bm{P}$ and $\bm{Q}$, $m \in \mathbb{N}$ and $a \in \mathcal{A}$,
define the quantity $H_{\bm{P},\bm{Q}}^{m,a} (\vartheta)$ by
\begin{equation*}
    H_{\bm{P},\bm{Q}}^{m,a} (\vartheta) = \sum_{v \in \vartheta^{m} (a)} - \mathbb{P} [\vartheta_{\bm{Q}}^m (a) = v] \log \mathbb{P} [\vartheta_{\bm{P}}^m (a) = v] \tp
\end{equation*}
Further, let $\bm{H}_{\bm{P},\bm{Q}}^m (\vartheta)$ denote the vector $(H_{\bm{P},\bm{Q}}^{m,a} (\vartheta))_{a \in \mathcal{A}}$.
We first prove some properties of the quantity $\bm{H}_{\bm{P},\bm{Q}}^{m} (\vartheta)$ which we will use in the proof of \cref{p:relative-local-dim}.

\begin{lemma}\label{l:fun-word-decomp}
    If $\vartheta$ is a primitive and compatible set-valued substitution and $\bm{P}$ and $\bm{Q}$ are permissible probabilities, then for all $m \in \mathbb{N}$, $a \in \mathcal{A}$ and $s \in \vartheta (a)$,
    \begin{equation*}
        \sum_{v \in \vartheta^m (s)} \mathbb{P} [\vartheta_{\bm{Q}}^m (s) = v] \log \mathbb{P} [\vartheta_{\bm{P}}^m (s) = v] = \sum_{b \in \mathcal{A}} \lvert \vartheta (a) \rvert_b \, H_{\bm{P},\bm{Q}}^{m,b} (\vartheta) \tp
    \end{equation*}
\end{lemma}
\begin{proof}
    Since $\vartheta$ is compatible, we can decompose each $v \in \vartheta^m (s)$ into inflation words $v = v^1 \cdots v^{\lvert \vartheta (a) \rvert}$.
    By the Markov property of $\vartheta_{\bm{P}}$ (respectively $\vartheta_{\bm{Q}}$),
    \begin{equation*}
        \mathbb{P} [\vartheta_{\bm{P}}^m (s) = v] = \mathbb{P} [\vartheta_{\bm{P}}^m(s_1) = v^1] \cdots \mathbb{P} [\vartheta_{\bm{P}}^m(s_{\lvert \vartheta (a) \rvert}) = v^{\lvert \vartheta (a) \rvert}] \text{.}
    \end{equation*}
    Therefore
    \begin{equation*}
        \begin{split}
            \sum_{v \in \vartheta^m (s)} \mathbb{P} [\vartheta_{\bm{Q}}^m (s) = v]&\log \mathbb{P} [\vartheta_{\bm{P}}^m (s) = v]=\\
                                                                                  &= \sum_{b \in \mathcal{A}} \lvert \vartheta (a) \rvert_b \sum_{w \in \vartheta^m (b)} \mathbb{P} [\vartheta_{\bm{Q}}^m (b) = w] \log \mathbb{P} [\vartheta_{\bm{P}}^m (b) = w] \\
                                                                                  &= \sum_{b \in \mathcal{A}} \lvert \vartheta (a) \rvert_b \, H_{\bm{P},\bm{Q}}^{m,b} (\vartheta)
                                                                                                                                  \tc
        \end{split}
    \end{equation*}
    which completes the proof.
\end{proof}

\begin{lemma}\label{l:fun-powers}
    If $\vartheta$ is a primitive and compatible set-valued substitution satisfying the disjoint set condition, with right Perron--Frobenius eigenvector $\bm{R}$, and $\bm{P}$ and $\bm{Q}$ are permissible probabilities, then
    \begin{equation*}
        \frac{1}{\lambda^m} \bm{H}_{\bm{P},\bm{Q}}^{m} (\vartheta) \cdot \bm{R} \rightarrow \frac{1}{\lambda-1} \bm{H}_{\bm{P},\bm{Q}}^{1} (\vartheta) \cdot \bm{R}
    \end{equation*}
    as $m \rightarrow \infty$.
\end{lemma}
\begin{proof}
    Since $\vartheta$ satisfies the disjoint set condition, for all $m \in \mathbb{N}$ and $a \in \mathcal{A}$,
    \begin{align*}
        \bm{H}_{\bm{P},\bm{Q}}^{m+1} (\vartheta) \cdot \bm{R}
        ={}& \sum_{a \in \mathcal{A}} R_a \sum_{v \in \vartheta^{m+1} (a)} \mathbb{P} [\vartheta_{\bm{Q}}^{m+1} (a) = v] \log \mathbb{P} [\vartheta_{\bm{P}}^{m+1} (a) = v] \\
        ={}& \sum_{a \in \mathcal{A}} R_a \sum_{s \in \vartheta (a)} \mathbb{P} [\vartheta_{\bm{Q}} (a) = s] \log \mathbb{P} [\vartheta_{\bm{P}} (a) = s]\\
           &+ \sum_{a \in \mathcal{A}} R_a \sum_{s \in \vartheta (a)} \mathbb{P} [\vartheta_{\bm{Q}} (a) = s] \sum_{v \in \vartheta^m (s)} \mathbb{P} [\vartheta_{\bm{Q}}^m (s) = v] \log \mathbb{P} [\vartheta_{\bm{P}}^m (s) = v] \\
        ={}& \bm{H}_{\bm{P}, \bm{Q}}^{1} (\vartheta) \cdot \bm{R} + \sum_{b \in \mathcal{A}} H_{\bm{P},\bm{Q}}^{m,b} (\vartheta) \sum_{a \in \mathcal{A}} \lvert \vartheta (a) \rvert_b R_a \\
        ={}& \bm{H}_{\bm{P}, \bm{Q}}^{1} (\vartheta) \cdot \bm{R} + \lambda \sum_{b \in \mathcal{A}} R_b H_{\bm{P},\bm{Q}}^{m,b} (\vartheta) \\
        ={}& \bm{H}_{\bm{P}, \bm{Q}}^{1} (\vartheta) \cdot \bm{R} + \lambda \bm{H}_{\bm{P}, \bm{Q}}^{m} (\vartheta) \cdot \bm{R}
        \tp
    \end{align*}
    In the second equality we use the Markov property of $\vartheta_{\bm{P}}$ and $\vartheta_{\bm{Q}}$, laws of logarithms, and that $\sum_{v \in \vartheta^m (s)} \mathbb{P} [\vartheta_{\bm{Q}}^m (s) = v] = 1$ for all $s \in \vartheta (a)$; in the third we apply \cref{l:fun-word-decomp}; in the fourth we use that $M_{\vartheta} \bm{R} = \lambda \bm{R}$.
    Applying the above inductively,
    \begin{equation*}
        \frac{1}{\lambda^m} \bm{H}_{\bm{P},\bm{Q}}^{m} (\vartheta) \cdot \bm{R} = \sum_{j=1}^{m} \frac{1}{\lambda^j} \bm{H}_{\bm{P},\bm{Q}}^{1} (\vartheta) \cdot \bm{R} \fto{m\to\infty} \frac{1}{\lambda-1} \bm{H}_{\bm{P},\bm{Q}}^{1} (\vartheta) \cdot \bm{R} \tc
    \end{equation*}
    which completes the proof.
\end{proof}

Any bi-infinite sequence $x$ in the subshift of a recognisable random substitution can be written as a bi-infinite concatenation of exact inflation words $(w^{n,a_n})$, where $w^{n,a_n}$ is an inflation word generated from the letter $a_n$.
Given a recognisable set-valued substitution $\vartheta$, $a \in \mathcal{A}$ and $w \in \vartheta (a)$, we define the inflation word frequency of $(a,w)$ in $x \in X_{\vartheta}$ by
\begin{align*}
    f_x(a,w) &= \lim_{n \to \infty} f^n_x(a,w)\\
    f^n_x(a,w) &= \frac{1}{2n+1} \# \{ m \colon a_m = a, w^{m,a_m} = w, w^{m,a_m} \mbox{ in } x_{[-n,n]}  \},
\end{align*}
provided the limit exists.
For a given frequency measure $\mu_{\bm{P}}$, the inflation word frequency of a $\mu_{\bm{P}}$-typical word is determined by the production probabilities.
More specifically, we have the following.
\begin{lemma}\label{l:expected-inf-frequencies}
    Let $\vartheta_{\bm{P}} = (\vartheta, \bm{P})$ be a primitive, compatible and recognisable random substitution with corresponding frequency measure $\mu_{\bm{P}}$.
    For $\mu_{\bm{P}}$-almost every $x \in X_{\vartheta}$ , the inflation word frequency exists and is given by
    \begin{equation*}
        f_x(a,w) = \frac{1}{\lambda} R_a \mathbb{P}[\vartheta_{\bm{P}} (a) = w],
    \end{equation*}
    for all $a \in \mathcal{A}$ and $w \in \vartheta(a)$.
\end{lemma}
\begin{proof}
    Let $A_{a,w}$ be the set of points $x \in X_{\vartheta}$ such that the above does \emph{not} hold.
    We show that $A_{a,w}$ is a null set.
    Taking the complement and then the intersection over all $a,w$ gives a full-measure set with the required property.
    Given $\varepsilon > 0$, let $E(n,\varepsilon)$ be the set of $x \in X_{\vartheta}$ such that
    \begin{equation*}
        \bigl|f^n_x(a,w) - \frac{1}{\lambda} R_a \mathbb{P}[\vartheta_{\bm{P}} (a) = w]\bigr| > \varepsilon.
    \end{equation*}
    By the Borel--Cantelli lemma, it suffices to show that
    \begin{equation*}
        \sum_{n \in \mathbb{N}} \mu_{\bm{P}}(E(n,\varepsilon)) < \infty
    \end{equation*}
    for all $\varepsilon > 0$ in order to conclude that $A_{a,w}$ is a nullset.
    To this end, we show that $\mu_{\bm{P}}(E(n,\varepsilon))$ decays exponentially with $n$.
    Given $u$ with $|u| = 2n+1 > 2 \kappa(\vartheta)$, let $u^R$ denote the recognisable core of $u$, which has length at least $\lvert u \rvert - 2 \kappa (\vartheta)$.
    \cref{l:recog-key-lemma} gives that
    \begin{equation*}
        \mu_{\bm{P}}([u]) \leq \frac{\kappa(\vartheta)}{\lambda} \mu_{\bm{P}}([v]) \mathbb{P}[\vartheta_{\bm{P}} (v) = u^R] \\
        = \frac{\kappa(\vartheta)}{\lambda} \mu_{\bm{P}}([v]) \prod_{i=1}^{|v|} \mathbb{P}[\vartheta_{\bm{P}} (v_i) = w^{i,v_i}]
    \end{equation*}
    where each $w^{i,v_i}$ is the inflated image of $v_i$ in $u^R$.
    By compatibility, we can choose an integer $N$ such that every $v$ of length at least $N$ satisfies $|v| (R_a - \varepsilon/3) \leq |v|_a \leq |v| (R_a + \varepsilon/3)$ for all $a \in \mathcal{A}$.
    For each $v$ and $a \in \mathcal{A}$, let $A_a (v)$ denote the set of $u' \in \vartheta (v)$ such that the frequency of indices $i \in \{ j: a_j = a \}$ with $w^{i,a} = w$ deviates from $\mathbb{P}[\vartheta_{\bm{P}}(a) = w]$ by more than $\varepsilon/3$.
    Since $\vartheta_{\bm{P}}$ acts independently on letters, it follows by Cramér's theorem that the sum
    $\sum_{u'\in A(v)} \mathbb{P}[\vartheta_{\bm{P}} (v) = u']$
    decays exponentially with $|v|_a$ (and hence with $|v|$).
    In particular, there is a constant $C > 0$, independent of the choice of $v$, such that
    \begin{equation}\label{eq:cramer}
        \sum_{u'\in A(v)} \mathbb{P}[\vartheta_{\bm{P}} (v) = u'] \leq e^{-C n} \text{.}
    \end{equation}
    Note that if $u$ is a sufficiently long legal word and has $[u] \cap E(n,\varepsilon) = \varnothing$, then we require that $u^R \in A(v)$.
    Indeed, if $u' \notin A(v)$ and $\lvert v \rvert \geq N$, then the relative inflation word frequency of $w$ is bounded above by
    \begin{equation*}
        \begin{split}
            \frac{\{ j: a_j = a \}}{|v|} \frac{|v|}{|u|}  \left( \mathbb{P}[\vartheta_{\bm{P}} (a) = w] + \frac{\varepsilon}{3} \right) &\leq \frac{1}{\lambda} \left( R_a + \frac{\varepsilon}{3} \right) \left( \mathbb{P}[\vartheta_{\bm{P}} (a) = w] + \frac{\varepsilon}{3} \right)\\
                                                                                                                                        &\leq \frac{1}{\lambda} R_a \mathbb{P}[\vartheta_{\bm{P}} (a) = w] + \varepsilon
        \end{split}
    \end{equation*}
    and, similarly, bounded below by $R_a \mathbb{P} [\vartheta_{\bm{P}} (a) = w]/\lambda - \varepsilon$; hence, $[u^R] \cap E(n,\varepsilon) = \varnothing$.
    Let $\mathcal{V}_n$ denote set of all words which appear as the (unique) preimage of the recognisable core of a word of length $n$.
    It then follows by \cref{l:recog-key-lemma} that
    \begin{equation*}
        \mu_{\bm{P}}(E(n,\varepsilon))
        \leqslant \sum_{\substack{u \in \mathcal{L}_{\vartheta}^{n} \\ [u] \cap E(n,\varepsilon) \neq \varnothing}} \mu_{\bm{P}}([u])
        \leq  \frac{\kappa(\vartheta)}{\lambda} \sum_{v \in \mathcal{V}_{n}} \mu_{\bm{P}}([v]) \sum_{u' \in A(v)} \mathbb{P}[\vartheta_{\bm{P}} (v) = u']
        \leqslant e^{-Cn},
    \end{equation*}
    where in the final inequality we have used \cref{eq:cramer} and that
    \begin{equation*}
        \sum_{v \in \mathcal{V}_n} \mu_{\bm{P}} ([v]) \leq \sum_{j=1}^n \sum_{v \in \mathcal{L}_{\vartheta}^{j}} \mu_{\bm{P}} ([v]) \leq n,
    \end{equation*}
    absorbing this contribution and the $\kappa(\vartheta)/\lambda$ factor into the constant $C$.
    It follows that
    \begin{equation*}
        \sum_{n=1}^\infty \mu_{\bm{P}}(E(n,\varepsilon)) \leq \sum_{n=1}^\infty e^{-Cn} < \infty \text{,}
    \end{equation*}
    and the result then follows by the Borel--Cantelli lemma.
\end{proof}

Finally, we require the following bounds on the exponential scaling rate of measures of cylinders, which is essentially a consequence of \cref{it:lq-inf-equiv}.
In particular, these give bounds on the possible local dimensions of the measure.
\begin{proposition}\label{p:meas-bounds}
    If $\vartheta_{\bm{P}}$ is a primitive and compatible random substitution with corresponding frequency measure $\mu_{\bm{P}}$, then there are values $0<s_1<s_2<\infty$ and $c_1,c_2>0$ such that for all $n\in\N$ and $v\in\mathcal{L}^n(X_\vartheta)=\mathcal{L}^n_\vartheta$,
    \begin{equation*}
        s_1\cdot n+c_1\leq\log\mu_{\bm{P}} ([v]) \leq s_2\cdot n+c_2
    \end{equation*}
\end{proposition}
\begin{proof}
    By \cref{it:lq-inf-equiv}, for all $k\in\N$ and $q>1$,
    \begin{equation*}
        \tau_{\mu_{\bm{P}}} (q) \leq \frac{1}{\lambda^k-1} \varphi_k(q)\text{;}
    \end{equation*}
    and for $q<0$,
    \begin{equation*}
        \frac{1}{\lambda^k-1}\varphi_k(q)\leq \tau_{\mu_{\bm{P}}}(q)\tc
    \end{equation*}
    Moreover, for each $k$, with
    \begin{align*}
        \beta_{k,\min}\coloneqq \lim_{q\to\infty}\frac{\varphi_k(q)}{q(\lambda^k-1)}&=-\frac{1}{\lambda^k-1}\sum_{a\in\mathcal{A}} R_a\log\left(\min_{v\in\vartheta^k(a)}\mathbb{P}[\vartheta_{\bm{P}}^k(a)=v]\right)\\
        \beta_{k,\max}\coloneqq \lim_{q\to-\infty}\frac{\varphi_k(q)}{q(\lambda^k-1)}&=-\frac{1}{\lambda^k-1}\sum_{a\in\mathcal{A}} R_a\log\left(\max_{v\in\vartheta^k(a)}\mathbb{P}[\vartheta_{\bm{P}}^k(a)=v]\right)\tc
    \end{align*}
    it follows that $[\beta_{k,\min},\beta_{k,\max}]\subset (0,\infty)$ is a decreasing nested sequence of intervals, so with $\beta_{\min}=\lim_{k\to\infty}\beta_{k,\min}$ and $\beta_{\max}=\lim_{k\to\infty}\beta_{k,\max}$,
    \begin{equation*}
        0<\beta_{\min}\leq\lim_{q\to\infty}\tau_{\mu_{\bm{P}}}(q)\leq\lim_{q\to -\infty}\tau_{\mu_{\bm{P}}}(q)\leq\beta_{\max}<\infty.
    \end{equation*}
    Applying \cref{l:lq-core}(b) gives the result.
\end{proof}
Finally, we obtain our main conclusion concerning relative local dimensions.
\begin{proposition}\label{p:relative-local-dim}
    Let $\vartheta$ be a primitive, compatible and recognisable set-valued substitution, let $\bm{P}$ and $\bm{Q}$ be permissible probabilities, and let $\mu_{\bm{P}}$ and $\mu_{\bm{Q}}$ denote the respective frequency measures.
    Then, for $\mu_{\bm{Q}}$-almost all $x \in X_{\vartheta}$,
    \begin{equation}\label{e:cross-ldim}
        \ldim(\mu_{\bm{P}}, x) = \frac{1}{\lambda-1} \sum_{a\in\mathcal{A}}R_a\sum_{v\in\vartheta(a)}-\mathbb{P} [\vartheta_{\bm{Q}}^m (a) = v] \log \mathbb{P} [\vartheta_{\bm{P}}^m (a) = v] \tp
    \end{equation}
\end{proposition}
\begin{proof}
    Fix $m\in\N$.
    It follows by \cref{l:recognisability-recursion} that since $\vartheta_{\bm{P}}$ is recognisable, so is $\vartheta_{\bm{P}}^m$.
    For each $x \in X_{\vartheta}$ and $n \in \mathbb{N}$ with $n > \kappa (\vartheta^m)$, let $u_{-}^{n} (x)$ denote the recognisable core of $x_{[-n,n]}$ and let $u_{+}^{n} (x)$ denote an inflation word of minimal length that contains $x_{[-n,n]}$.
    By compatibility, $\lvert u_{-}^{n} (x) \rvert / (2n+1) \rightarrow \lambda^{-m}$ and $\lvert u_{+}^{n} (x) \rvert / (2n+1) \rightarrow \lambda^{-m}$ as $n \rightarrow \infty$.
    Further, let $v_{-}^{n} (x)$ be the legal word such that $u_{-}^{n} (x) \in \vartheta^m (v_{-}^{n} (x))$ and $v_{+}^{n} (x)$ be the legal word such that $u_{+}^{n} (x) \in \vartheta^m (v_{+}^{n} (x))$.
    Then, it follows by \cref{l:recog-key-lemma} and the definition of local dimension that
    \begin{align*}
        \liminf_{n \rightarrow \infty} &\left( -\frac{1}{2n+1} \log \mu_{\bm{P}} ([u_{-}^{n} (x)]) - \frac{1}{2n+1} \log \mathbb{P} [\vartheta_{\bm{P}} (v_{-}^{n} (x)) = u_{-}^{n} (x)] \right)\\
                                       &\leq \lldim (\mu_{\bm{P}}, x) \leq \uldim (\mu_{\bm{P}},x)\\
                                       &\leq \limsup_{n \rightarrow \infty} \left( -\frac{1}{2n+1} \log \mu_{\bm{P}} ([u_{+}^{n} (x)]) - \frac{1}{2n+1} \log \mathbb{P} [\vartheta_{\bm{P}} (v_{+}^{n} (x)) = u_{+}^{n} (x)] \right) \tp
    \end{align*}
    By \cref{p:meas-bounds}, there exists a constant $C \geq 0$ such that for all $x \in X_{\vartheta}$,
    \begin{equation*}
        0 \leq \liminf_{n \rightarrow \infty} -\frac{1}{2n+1} \log \mu_{\bm{P}} ([u_{-}^{n} (x)]) \leq \limsup_{n \rightarrow \infty} -\frac{1}{2n+1} \log \mu_{\bm{P}} ([u_{+}^{n} (x)]) \leq C \tp
    \end{equation*}
    Hence, it follows from the above that
    \begin{equation}\label{eq:loc-dim-bounds}
        \begin{aligned}
            \liminf_{n \rightarrow \infty} - &\frac{1}{2n+1} \log \mathbb{P} [\vartheta_{\bm{P}} (v_{-}^{n} (x)) = u_{-}^{n} (x)]\\
                                             &\leq \lldim(\mu_{\bm{P}}, x)\leq \uldim(\mu_{\bm{P}},x)\\
                                             &\leq \limsup_{n \rightarrow \infty} - \frac{1}{2n+1} \log \mathbb{P} [\vartheta_{\bm{P}} (v_{+}^n (x)) = u_{+}^{n} (x)] + \frac{C}{\lambda^m} \tp
        \end{aligned}
    \end{equation}

    We now show that for $\mu_{\bm{Q}}$-almost all $x \in X_{\vartheta}$,
    \begin{align*}
        \liminf_{n \rightarrow \infty} - \frac{1}{n} \log \mathbb{P} [\vartheta_{\bm{P}} (v_{-}^{n} (x)) = u_{-}^{n} (x)] &= \limsup_{n \rightarrow \infty} - \frac{1}{n} \log \mathbb{P} [\vartheta_{\bm{P}} (v_{+}^n (x)) = u_{+}^{n} (x)]\\
                                                                                                                          &= \frac{1}{\lambda^m} \bm{H}_{\bm{P},\bm{Q}}^{m} (\vartheta) \cdot \bm{R} \tp
    \end{align*}
    By compatibility, we can decompose the production probabilities into inflation tiles as
    \begin{equation*}
        \mathbb{P} [\vartheta_{\bm{P}}^m (v_{-}^{n} (x)) = u_{-}^{n} (x)] = \prod_{a \in \mathcal{A}} \prod_{w \in \vartheta^{m} (a)} \mathbb{P} [\vartheta_{\bm{P}}^m (a) = w]^{N_{a,w} (x,n)} \tc
    \end{equation*}
    where, for each $a \in \mathcal{A}$ and $w \in \vartheta^m (a)$, $N_{a,w} (x,n)$ denotes the number of $a$'s in $v_{-}^{n} (x)$ which map to $w$.
    It follows by \cref{l:expected-inf-frequencies}, applied to $\vartheta_{\bm{Q}}^m$, that for $\mu_{\bm{Q}}$-almost all $x \in X_{\vartheta}$,
    \begin{equation*}
        \frac{1}{2n+1} N_{a,w} (x,n) \rightarrow \frac{1}{\lambda^m} R_a \mathbb{P} [\vartheta_{\bm{Q}}^m (a) = w]
    \end{equation*}
    for all $a \in \mathcal{A}$ and $w \in \vartheta^m (a)$.
    Hence, it follows that
    \begin{align*}
        \lim_{n\to\infty}- \frac{1}{2n+1} \log \mathbb{P}&[\vartheta_{\bm{P}}^m (v_{-}^{n} (x)) = u_{-}^{n} (x)]\\
                                                         &= \frac{1}{\lambda^m} \sum_{a \in \mathcal{A}} R_a \sum_{v \in \vartheta^m (a)} \mathbb{P} [\vartheta_{\bm{Q}}^m (a) = v] \log \mathbb{P} [\vartheta_{\bm{P}}^m (a) = v]\\
                                                         &= \frac{1}{\lambda^m} \bm{H}_{\bm{P}, \bm{Q}}^m (\vartheta) \cdot \bm{R} \tc
    \end{align*}
    with the same convergence holding for $u_{+}^{n} (x)$ by identical arguments.
    Thus, it follows from \cref{eq:loc-dim-bounds} that
    \begin{equation*}
        \frac{1}{\lambda^m} \bm{H}_{\bm{P},\bm{Q}}^m (\vartheta) \cdot \bm{R} \leq \lldim (\mu_{\bm{P}}, x) \leq \uldim (\mu_{\bm{P}}, x) \leq \frac{1}{\lambda^m} \bm{H}_{\bm{P},\bm{Q}}^m (\vartheta) \cdot \bm{R} + \frac{C}{\lambda^m} \tp
    \end{equation*}
    Since the above holds for all $m \in \mathbb{N}$, by letting $m \rightarrow \infty$ it follows by \cref{l:fun-powers} that $\ldim(\mu_{\bm{P}}, x)$ exists and
    \begin{equation*}
        \ldim (\mu_{\bm{P}}, x) = \frac{1}{\lambda-1} \bm{H}_{\bm{P},\bm{Q}}^1 (\vartheta) \cdot \bm{R} \tc
    \end{equation*}
    which completes the proof.
\end{proof}

\subsection{Proof of the multifractal formalism}\label{ss:multi-proof}
In this section, we apply the results obtained in the previous section, along with results on the $L^q$-spectrum under recognisability, to prove \cref{it:multi-formalism}.
\begin{proof}[of \cref{it:multi-formalism}]
    We first obtain the results for the $L^q$-spectrum.
    Since every recognisable random substitution satisfies the disjoint set condition, \cref{p:inflation-isc-dsc} gives that $T_{\vartheta,\bm{P}} (q) = (\lambda-1)^{-1} \varphi_1 (q)$ for all $q \in \mathbb{R}$.
    If $q<0$, then by \cref{it:lq-inf-equiv} and \cref{p:q<0-lower-recognisable},
    \begin{equation*}
        \frac{1}{\lambda-1} \varphi_1 (q) = T_{\vartheta,\bm{P}} (q) \leq \tau_{\mu_{\bm{P}}} (q) \leq \overline{\tau}_{\mu_{\bm{P}}} (q) \leq \frac{1}{\lambda-1} \varphi_1 (q) \tc
    \end{equation*}
    so we conclude that $\tau_{\mu_{\bm{P}}} (q)$ exists and equals $(\lambda-1)^{-1} \varphi_1 (q)$.
    For $q \geq 0$, the result follows already from \cref{ic:lq-lim-eq}.

    We now obtain the results on the multifractal spectrum.
    In light of \cref{p:multi-upper}, it remains to show that $f_{\mu_{\bm{P}}}(\alpha)\geq\tau_{\mu_{\bm{P}}}^*(\alpha)$ for each $\alpha\in\R$.
    As proved above, for all $q\in\R$,
    \begin{equation*}
        \tau_{\mu_{\bm{P}}}(q) = \frac{1}{\lambda - 1} \varphi_1 (q)= \frac{1}{\lambda-1}\sum_{a \in \mathcal{A}} R_a T_a(q)
    \end{equation*}
    where for each $a\in\mathcal{A}$
    \begin{equation*}
        T_a(q)=-\log\sum_{s\in\vartheta(a)}\mathbb{P}[\vartheta_{\bm{P}}(a)=s]^q\tp
    \end{equation*}

    First, fix $\alpha\in(\alpha_{\min},\alpha_{\max})$ and let $q\in\R$ be chosen so that $\tau_{\mu_{\bm{P}}}'(q)=\alpha$.
    Observe that $q\alpha-\tau_{\mu_{\bm{P}}}(q)=\tau_{\mu_{\bm{P}}}^*(\alpha)$.
    Then define $\bm{Q}$ by the rule
    \begin{equation*}
        \mathbb{P}[\vartheta_{\bm{Q}}(a)=s]= \mathbb{P}[\vartheta_{\bm{P}}(a)=s]^q e^{T_a(q)}
    \end{equation*}
    for all $a\in\mathcal{A}$ and $s\in\vartheta(a)$.
    Then by \cref{ic:entropy-recover},
    \begin{align*}
        \dimH\mu_{\bm{Q}} ={}& \frac{1}{\lambda-1}\sum_{a\in\mathcal{A}}R_a\Bigl(-\sum_{v \in \vartheta (a)} \mathbb{P} [\vartheta_{\bm{Q}} (a) = v] \log \mathbb{P} [\vartheta_{\bm{Q}} (a) = v]\Bigr)\\
        ={}&q\cdot\frac{1}{\lambda-1}\sum_{a\in\mathcal{A}}R_a\Bigl(-\sum_{v \in \vartheta (a)} \mathbb{P} [\vartheta_{\bm{Q}} (a) = v] \log \mathbb{P} [\vartheta_{\bm{P}}(a) = v]\Bigr)\\
                              &-\frac{1}{\lambda-1}\sum_{a\in\mathcal{A}}R_aT_a(q)\sum_{v \in \vartheta (a)} \mathbb{P} [\vartheta_{\bm{Q}} (a) = v]\\
                              ={}&q\alpha-\tau_{\mu_{\bm{P}}}(q)=\tau_{\mu_{\bm{P}}}^*(\alpha)
    \end{align*}
    since
    \begin{align*}
        \tau_{\mu_{\bm{P}}}'(q)&= \frac{1}{\lambda-1}\sum_{a\in\mathcal{A}}R_a\frac{-\sum_{v \in \vartheta (a)} \mathbb{P} [\vartheta_{\bm{P}} (a) = v]^q \log \mathbb{P} [\vartheta_{\bm{P}}(a) = v]}{e^{-T_a(q)}}\\
                    &=\frac{1}{\lambda-1}\sum_{a\in\mathcal{A}}R_a\Bigl(-\sum_{v \in \vartheta (a)} \mathbb{P} [\vartheta_{\bm{Q}} (a) = v] \log \mathbb{P} [\vartheta_{\bm{P}}(a) = v]\Bigr)\tp
    \end{align*}
    In fact, this shows that $\ldim(\mu_{\bm{P}},x)=\alpha$ for $\mu_{\bm{Q}}$-almost all $x\in X_\vartheta$ by \cref{p:relative-local-dim}.
    Thus $f_{\mu_{\bm{P}}}(\alpha)\geq\dimH\mu_{\bm{Q}}=\tau_{\mu_{\bm{P}}}^*(\alpha)$, as required.

    The result for $\alpha=\alpha_{\min}$ (resp. $\alpha=\alpha_{\max}$) follows similarly by taking a degenerate probability vector $\bm{Q}$ assigning equal value to the realisations of $\vartheta(a)$ with maximal (resp. minimal) probabilities given by $\bm{P}$, and zero otherwise.
    The corresponding non-degenerate sub-substitution is also compatible and recognisable, so the same arguments yield the corresponding bounds.
\end{proof}

\section{Examples, counterexamples and applications}\label{s:examples}
\subsection{Failure of bounds for negative \texorpdfstring{$q$}{q} without recognisability}
In the following two examples, we show the results in \cref{it:lq-inf-equiv} do not extend in general to give an upper bound for the $L^q$-spectrum in terms of the inflation word $L^q$-spectrum, for $q<0$.
In \cref{ex:q<0-failure-full-shift}, we construct a class of frequency measures on the full-shift on two letters for which the $L^q$-spectrum and inflation word analogue differ in the $q<0$ case.
The random substitutions that give rise to these frequency measures are not compatible, but in \cref{ex:q<0-non-inflation} we present a compatible analogue.

In contrast, in \cref{ex:isc-q<0-equal}, we give an example showing that the results for $q<0$ can hold for all $q\in\R$ under the identical set condition with identical production probabilities.
\begin{example}\label{ex:q<0-failure-full-shift}
    Let $p_1 < p_2 \in (0,1)$ such that $p_1 + 3 p_2 = 1$ and let $\vartheta_{\bm{P}}$ be the random substitution defined by
    \begin{equation*}
        \vartheta_{\bm{P}} \colon a, b \mapsto
        \begin{cases}
            ab &\text{ with probability $p_1$}\\
            ba &\text{ with probability $p_2$}\\
            aa &\text{ with probability $p_2$}\\
            bb &\text{ with probability $p_2$}\\
        \end{cases}
    \end{equation*}
    We show for all sufficiently small $q < 0$ that $\tau_{\mu_{\bm{P}}} (q) > T_{\vartheta,\bm{P}}(q)$.
    Observe that, for each $k \in \mathbb{N}$, the word $v^k = (ab)^{2^k} \in \vartheta^{k+1} (a) \cap \vartheta^{k+1} (b)$ occurs with probability
    \begin{equation*}
        \mathbb{P} [\vartheta_{\bm{P}}^{k+1} (a) = v^k] = \mathbb{P} [\vartheta_{\bm{P}}^{k+1} (b) = v^k] = p_1^{2^{k}} \tp
    \end{equation*}
    Clearly, this is the minimal possible probability with which a level-$k$ inflation word can occur, so it follows that
    \begin{equation*}
        \lim_{q \rightarrow - \infty} \frac{T_{\vartheta,\bm{P}}(q)}{q} = -\frac{1}{2} \log p_1 \tp
    \end{equation*}
    Now, let $u \in \mathcal{L}_{\vartheta}^{2^{k+1}}$ be arbitrary.
    We show that $\mu_{\bm{P}} ([u]) \geq p_1^{2^{k-1}} p_2^{2^{k-1}} / 2$.
    Since $\vartheta (a) = \vartheta (b)$ with identical production probabilities, it follows by \cref{l:key-lemma} that for any choice of $w \in \mathcal{L}_{\vartheta}^{2^k + 1}$
    \begin{equation*}
        \mu_{\bm{P}} ([u]) = \frac{1}{2} \left( \mathbb{P} [\vartheta_{\bm{P}} (w)_{[1,2^{k+1}]} = u] + [\vartheta_{\bm{P}} (w)_{[2,2^{k+1}+1]} = u]  \right) \tp
    \end{equation*}
    If $\mathbb{P} [\vartheta_{\bm{P}} (w)_{[1,2^{k+1}]} = u] \geq p_1^{2^{k-1}} p_2^{2^{k-1}}$, then we are done, otherwise at least half of the letters in $v$ must be sent to $ab$.
    But then for $u$ to appear from the second letter, at least half of the letters in $v$ must be sent to $ba$ or $bb$, so $\mathbb{P} [\vartheta_{\bm{P}} (w)_{[2,2^{k+1} + 1]} = u] \geq p_1^{2^{k-1}} p_2^{2^{k-1}}$.
    Hence, $\mu_{\bm{P}} ([u]) \geq p_1^{2^{k-1}} p_2^{2^{k-1}} / 2$ so, in particular,
    \begin{equation*}
        \min_{u \in \mathcal{L}_{\vartheta}^{2^{k+1}}} \mu_{\bm{P}} ([u]) \geq \frac{1}{2} p_1^{2^{k-1}} p_2^{2^{k-1}} \tp
    \end{equation*}
    It follows that
    \begin{equation*}
        \lim_{q \rightarrow - \infty} \frac{\tau_{\mu_{\bm{P}}} (q)}{q} \leq - \frac{1}{4} (\log p_1 + \log p_2) < - \frac{1}{2} \log p_1 = \lim_{q \rightarrow - \infty} \frac{T_{\vartheta,\bm{P}}(q)}{q} \tp
    \end{equation*}
\end{example}

By a slight modification of this example, we can construct a compatible random substitution for which the two notions do not coincide.
\begin{example}\label{ex:q<0-non-inflation}
    Let $p_1 < p_2 \in (0,1)$ such that $p_1 + 3 p_2 = 1$ and let $\vartheta_{\bm{P}}$ be the random substitution defined by
    \begin{equation*}
        \vartheta_{\bm{P}} \colon a, b \mapsto
        \begin{cases}
            ab \, ba &\text{with probability }p_1\\
            ba \, ab &\text{with probability }p_2\\
            ab \, ab &\text{with probability }p_2\\
            ba \, ba &\text{with probability }p_2
        \end{cases}
    \end{equation*}
    By similar arguments to the previous example, we obtain
    \begin{equation*}
        \lim_{q \rightarrow - \infty} \frac{\tau_{\mu_{\bm{P}}} (q)}{q} \leq - \frac{1}{8} (\log p_1 + \log p_2) < - \frac{1}{4} \log p_1 = \lim_{q \rightarrow - \infty} \frac{T_{\vartheta,\bm{P}}(q)}{q} \tp
    \end{equation*}

\end{example}

The random substitution in \cref{ex:q<0-non-inflation} satisfies the identical set condition with identical production probabilities.
These conditions are also satisfied by the following example.
However, here the $L^q$-spectrum and inflation word analogue coincide for all $q\in\R$ by a direct argument.
\begin{example}\label{ex:isc-q<0-equal}
    We show that for the random substitution
    \begin{equation*}
        \vartheta_{\bm{P}} \colon a, b \mapsto
        \begin{cases}
            ab & \text{with probability }p\\
            ba & \text{with probability }1-p
        \end{cases}
    \end{equation*}
    the limit defining $\tau_{\mu_{\bm{P}}}(q)$ exists for all $q \in\R$, and
    \begin{equation*}
        \tau_{\mu_{\bm{P}}} (q) = T_{\vartheta,\bm{P}}(q) = \frac{1}{\lambda} \varphi_1 (q) = - \frac{1}{2} \log(p^q + (1-p)^q) \tp
    \end{equation*}
    \cref{ic:lq-lim-eq} gives the result for all $q > 0$ and that $\tau_{\mu_{\bm{P}}}(q) \geq T_{\vartheta,\bm{P}}(q) = 2^{-1} \varphi_1 (q)$ for all $q < 0$, so it only remains to verify for all $q<0$ that
    \begin{equation*}
        \overline{\tau}_{\mu_{\bm{P}}} (q) \leq T_{\vartheta,\bm{P}} (q).
    \end{equation*}

    Since $\vartheta (v^1) = \vartheta (v^2)$ for all $v^1, v^2 \in \mathcal{L}_{\vartheta}$, it follows from \cref{l:key-lemma} that for all $u \in \mathcal{L}_{\vartheta}^{2m}$ and any $v \in \mathcal{L}_{\vartheta}^{m+1}$,
    \begin{equation*}
        \mu_{\bm{P}} ([u]) = \frac{1}{2} \left(\mathbb{P} [\vartheta (v)_{[1,1+2 m-1]} =u] +\mathbb{P} [\vartheta (v)_{[2,2+2 m-1]} =u]\right)\tp
    \end{equation*}
    Let $\mathcal{V}_{2 m} = \{ (ab)^m, (ba)^m \}$.
    If $u \in \mathcal{L}_{\vartheta}^{2 m} \setminus \mathcal{V}_{2 m}$, then $u$ must contain $bb$ as a subword.
    This uniquely determines the cutting points in any inflation word decomposition, so there exists a unique $v$ and $j(u) \in \{1, 2\}$ such that $u \in \vartheta (v)_{[j(u), 2m + j(u) - 1]}$.
    It follows that
    \begin{equation*}
        \begin{split}
            \sum_{u \in \mathcal{L}_{\vartheta}^{2 m}} \mu_{\bm{P}} ([u])^q &\geq \sum_{u \in \mathcal{L}_{\vartheta}^{2 m} \setminus \mathcal{V}_{2 m}} \left( \frac{1}{2} \mathbb{P} [\vartheta_{\bm{P}} (v)_{[j(u), j(u) + 2 m - 1]}=u] \right)^q\\
                                                                   &\geq \frac{1}{2^q} \sum_{u \in \mathcal{L}_{\vartheta}^{2 m} \setminus \mathcal{V}_{2 m}} \mathbb{P} [\vartheta_{\bm{P}} (v_2 \cdots v_m) = u_{[3-j(u),2-j(u)+2m]} ]^q \tp
        \end{split}
    \end{equation*}
    Now, for every $w \in \vartheta (v_2 \cdots v_m)$ there is a $u$ such that $w = u_{[3-j(u),2-j(u)+2m]}$.
    Hence,
    \begin{equation*}
        \sum_{u \in \mathcal{L}^{2 m}_\vartheta} \mu_{\bm{P}} ([u])^q \geq \frac{1}{2^q} \sum_{w \in \vartheta (v_2 \cdots v_m)} \mathbb{P} [\vartheta_{\bm{P}} (v_2 \cdots v_m) = w]^q
    \end{equation*}
    and the conclusion follows by similar arguments to those used in the proofs of the main theorems.
\end{example}

\subsection{Examples with recognisability}
We first provide examples of random substitutions for which the multifractal formalism holds.
\begin{example}\label{ex:recog}
    Let $p>0$ and let $\vartheta_p$ be the random substitution defined by
    \begin{equation*}
        \vartheta_p \colon
        \begin{cases}
            a \mapsto
            \begin{cases}
                abb &\text{with probability $p$}\\
                bab &\text{with probability $1-p$}
            \end{cases}\\
            b \mapsto aa
        \end{cases}
    \end{equation*}
    Certainly $\vartheta_p$ is compatible, with corresponding primitive substitution matrix
    \begin{equation*}
        M=\begin{pmatrix}1&2\\2&0\end{pmatrix}\tc
    \end{equation*}
    Perron--Frobenius eigenvalue $(1+\sqrt{17})/2$, and right Perron--Frobenius eigenvector
    \begin{equation*}
        \left(\frac{-3+\sqrt{17}}{2},\frac{5-\sqrt{17}}{2}\right)\tp
    \end{equation*}
    One can verify that $\vartheta$ is recognisable since every occurrence of $aa$ intersects an image of $b$ and the adjacent letters then determine the cutting points.
    Thus by \cref{it:multi-formalism}, for all $q \in \mathbb{R}$
    \begin{equation*}
        \tau_{\mu_p} (q) = T_{\vartheta,\bm{P}}(q) = \frac{1}{\lambda-1} \varphi_1 (q) = -\frac{7-\sqrt{17}}{8}\log(p^q + (1-p)^q)
    \end{equation*}
    and measure $\mu_p$ satisfies the multifractal formalism.
    The asymptotes have slopes $-(7-\sqrt{17})\log(p)/8$ and $-(7-\sqrt{17})\log(1-p)/8$.
    A plot of the $L^q$-spectra and multifractal spectra for two choices of $p$ is given in \cref{f:recog}.

    For $p=1/2$, the $L^q$-spectrum of the measure $\mu_p$ is a straight line and the multifractal spectrum is equal to $\htop(X_{\vartheta})$ at $\htop(X_{\vartheta})$, and $-\infty$ otherwise.
\end{example}
\begin{figure}[t]
    \centering
    \begin{subcaptionblock}{.47\textwidth}
        \centering
        \begin{tikzpicture}[>=stealth,xscale=0.5,yscale=0.8,font=\tiny]
    \begin{scope}[thick,gray,->]
        \draw (-6.1,0) -- (6.1,0);
        \draw (0,-3.3) -- (0,1);
    \end{scope}

    \matrix [draw,above left,scale=0.8,row sep=0.3em, column sep=0.3em, inner sep=0.3em, align=left] at (current bounding box.south east) {
        \draw[thick] (0,0) -- (0.7,0); & \node[anchor=west,inner sep=0pt] {$\tau_{1/5}$}; \\
        \draw[thick, dashed] (0,0) -- (0.7,0); & \node[anchor=west,inner sep=0pt] {$\tau_{2/5}$}; \\
    };

    \begin{scope}[thick, smooth, variable=\x,domain=-6:6]
        \draw plot ({\x}, {-0.3596118*ln(0.2^\x+(1-0.2)^\x)});
        \draw[dashed] plot ({\x}, {-0.3596118*ln(0.4^\x+(1-0.4)^\x)});
    \end{scope}
\end{tikzpicture}
        \caption{$L^q$-spectra}
    \end{subcaptionblock}%
    \begin{subcaptionblock}{.47\textwidth}
        \centering
        \begin{tikzpicture}[>=stealth,xscale=10,yscale=12,font=\tiny]
    \begin{scope}[thick]
        \draw[gray, ->] (-0.05,0) -- (0.6,0);
        \draw[gray, ->] (0,-0.0416) -- (0,0.25);
        \draw plot file {figures/recog_1_5_pts.txt};
        \draw[dashed] plot file {figures/recog_2_5_pts.txt};
    \end{scope}

    \matrix [draw,fill=white,below left,scale=0.8,row sep=0.3em, column sep=0.3em, inner sep=0.4em, align=left] at (current bounding box.north east) {
        \draw[thick] (0,0) -- (0.7,0); & \node[anchor=west,inner sep=0pt] {$\tau_{1/5}^*$}; \\
        \draw[thick, dashed] (0,0) -- (0.7,0); & \node[anchor=west,inner sep=0pt] {$\tau_{2/5}^*$}; \\
    };
\end{tikzpicture}
        \caption{Multifractal spectra}
    \end{subcaptionblock}%
    \caption{$L^q$-spectra and multifractal spectra corresponding to a recognisable substitution for $p\in\{1/5,2/5\}$.}\label{f:recog}
\end{figure}
In the following example, we highlight that the multifractal spectrum need not have value 0 at the endpoints.
\begin{example}
    Let $\vartheta_{p}$ be the random substitution defined by
    \begin{equation*}
        \vartheta_p \colon
        \begin{cases}
            a \mapsto
            \begin{cases}
                abb &\text{with probability $p$}\\
                bab &\text{with probability $p$}\\
                bba &\text{with probability $1-2p$}
            \end{cases}\\
            b \mapsto aaa
        \end{cases}
    \end{equation*}
    Similarly to \cref{ex:recog}, $\vartheta_{p}$ is primitive, compatible and recognisable.
    Hence, \cref{it:multi-formalism} gives that
    \begin{equation*}
        \tau_{\mu_p} (q) = -\frac{3}{10} \log( 2p^q + (1-2p)^q) \tp
    \end{equation*}
    The asymptotes have slopes $-3\log(p)/10$ and $-3\log(1-2p)/10$.
    For $p=1/5$ and $p=2/5$, the $L^q$-spectrum and multifractal spectrum of $\mu_p$ are plotted in \cref{f:nonzero-example}.
    Here, we highlight that the endpoints of the multifractal spectrum need not be equal to zero.
\end{example}
\begin{figure}[t]
    \centering
    \begin{subcaptionblock}{.47\textwidth}
        \centering
        \begin{tikzpicture}[>=stealth,xscale=0.6,yscale=0.9,font=\tiny]
    \begin{scope}[thick,gray,->]
        \draw (-5.1,0) -- (5.1,0);
        \draw (0,-3) -- (0,1);
    \end{scope}

    \matrix [draw,above left,scale=0.8,row sep=0.3em, column sep=0.3em, inner sep=0.3em, align=left] at (current bounding box.south east) {
        \draw[thick] (0,0) -- (0.7,0); & \node[anchor=west,inner sep=0pt] {$\tau_{1/5}$}; \\
        \draw[thick, dashed] (0,0) -- (0.7,0); & \node[anchor=west,inner sep=0pt] {$\tau_{2/5}$}; \\
        };

    \begin{scope}[thick, smooth, variable=\x,domain=-5:5]
        \draw plot ({\x}, {-0.3*ln(2*(0.2^\x)+(1-0.4)^\x)});
        \draw[dashed] plot ({\x}, {-0.3*ln(2*(0.4^\x)+(1-0.8)^\x)});
    \end{scope}
\end{tikzpicture}
        \caption{$L^q$-spectra}
    \end{subcaptionblock}%
    \begin{subcaptionblock}{.47\textwidth}
        \centering
        \begin{tikzpicture}[>=stealth,xscale=10,yscale=8,font=\tiny]
    \begin{scope}[thick]
        \draw[gray, ->] (-0.04,0) -- (0.6,0);
        \draw[gray, ->] (0,-0.05) -- (0,0.35);
        \draw plot file {figures/nonzero_1_5_pts.txt};
        \draw[dashed] plot file {figures/nonzero_2_5_pts.txt};
    \end{scope}

    \matrix [fill=white,draw,below right,scale=0.8,row sep=0.3em, column sep=0.3em, inner sep=0.4em, align=left] at (current bounding box.north west) {
        \draw[thick] (0,0) -- (0.7,0); & \node[anchor=west,inner sep=0pt] {$\tau_{1/5}^*$}; \\
        \draw[thick, dashed] (0,0) -- (0.7,0); & \node[anchor=west,inner sep=0pt] {$\tau_{2/5}^*$}; \\
        };
\end{tikzpicture}
        \caption{Multifractal spectra}
    \end{subcaptionblock}%
    \caption{$L^q$-spectra and multifractal spectra corresponding to a recognisable substitution for $p\in\{1/5,2/5\}$.}\label{f:nonzero-example}
\end{figure}
\begin{example}\label{ex:intrinsic}
    Consider the random substitution on the three-letter alphabet $\mathcal{A}=\{a,b,c\}$ defined by
    \begin{equation*}
        \vartheta_{\mathbf{P}} \colon
        \begin{cases}
            a \mapsto
            \begin{cases}
                bbc &\text{with probability $p_1$}\\
                cbb &\text{with probability $1-p_1$}
            \end{cases}\\
            b \mapsto
            \begin{cases}
                cca &\text{with probability $p_2$}\\
                acc &\text{with probability $1-p_2$}
            \end{cases}\\
            c \mapsto
            \begin{cases}
                aab &\text{with probability $p_3$}\\
                baa &\text{with probability $1-p_3$}
            \end{cases}
        \end{cases}
    \end{equation*}
    for $p_1$, $p_2$, and $p_3$ in $(0,1)$.
    It is immediate that this substitution is compatible, and by considering the occurrences of 2, 3, or 4 letter repetitions, we see that this substitution is also recognisable.
    Moreover, the hypotheses of \cite[Theorem~4.8]{gmrs2023} are satisfied since $\vartheta$ is constant length and $\#\vartheta(a)=\#\vartheta(b)=\#\vartheta(c)$.
    In particular, the corresponding subshift $X_\vartheta$ is intrinsically ergodic with unique measure of maximal entropy given by taking $p_1=p_2=p_3=1/2$.

    It follows from \cite[Lemma~4.12]{gmrs2023} that the measure of maximal entropy is not a Gibbs measure with respect to the zero potential, so the system does not satisfy the usual specification property.
    For this choice of uniform probabilities, the $L^q$-spectrum is a straight line passing through the point $(1,0)$ with slope $\htop(X_\vartheta)=\log(2)/2$.
    More generally, the $L^q$-spectrum is given for all $q\in\R$ by the formula
    \begin{equation*}
        \tau_{\mu_{\bm{P}}}(q)=-\frac{1}{6}\left(\log\bigl((1 - p_1)^q + p_1^q\bigr)+ \log\bigl((1 - p_2)^q + p_2^q\bigr)+ \log\bigl((1 - p_3)^q + p_3^q\bigr)\right)
    \end{equation*}
    and the multifractal formalism is satisfied.

    For an example on an alphabet of size two, one may consider the random substitution
    \begin{equation*}
        \vartheta_{\mathbf{P}} \colon
        \begin{cases}
            a \mapsto
            \begin{cases}
                ababbb &\text{with probability $p_1$}\\
                abbabb &\text{with probability $1-p_1$}
            \end{cases}\\
            b \mapsto
            \begin{cases}
                baabaa &\text{with probability $p_2$}\\
                babaaa &\text{with probability $1-p_2$}
            \end{cases}
        \end{cases}
    \end{equation*}
    for $p_1$ and $p_2$ in $(0,1)$.
    The analysis of this example proceeds identically as above.
\end{example}

\subsection{Examples without recognisability}
Finally, we consider the two most commonly studied examples of random substitutions: random period doubling and random Fibonacci.
\begin{example}\label{ex:rpd-end}
    Given $p \in (0,1)$, let $\vartheta_p$ be the random period doubling substitution defined by
    \begin{equation*}
        \vartheta_p \colon
        \begin{cases}
            a \mapsto
            \begin{cases}
                ab &\text{with probability } p\\
                ba &\text{with probability } 1-p
            \end{cases}\\
            b \mapsto aa
        \end{cases}
    \end{equation*}
    and let $\mu_p$ denote the corresponding frequency measure.
    The substitution $\vartheta_p$ satisfies the disjoint set condition, so for all $q \in [0, \infty)$,
    \begin{equation*}
        \tau_{\mu_p} (q) = - \frac{2}{3} \log (p^q + (1-p)^q) \tp
    \end{equation*}
    The asymptote as $q\to\infty$ has slope $-2\log(\max\{p, 1-p\})/3$, which gives a sharp lower bound on the local dimensions of $\mu_p$.

    If $p = 1/2$, then the measure $\mu_p$ has linear $L^q$-spectrum for $q\geq 0$ given by
    \begin{equation*}
        \tau_{\mu_{1/2}} (q) = \frac{2}{3} (q-1) \log 2.
    \end{equation*}
    Since the substitution satisfies the disjoint set condition but is not recognisable, our results do not give the $L^q$-spectrum for $q<0$.
\end{example}

\begin{example}\label{ex:fib>0}
    The random Fibonacci substitution $\vartheta_p$ defined by
    \begin{equation*}
        \vartheta_{p} \colon
        \begin{cases}
            a \mapsto
            \begin{cases}
                ab &\text{with probability }p\\
                ba &\text{with probability }1-p
            \end{cases}\\
            b \mapsto a
        \end{cases}
    \end{equation*}
    does not satisfy either the identical set condition nor the disjoint set condition.
    Hence, we cannot apply \cref{ic:lq-lim-eq} to obtain a closed-form formula for $\tau_{\mu_p} (q)$.
    However, we can still apply \cref{it:lq-inf-equiv} to obtain a sequence of lower and upper bounds.
    The case $k=1$ gives the following bounds for all $0<q<1$:
    \begin{equation*}
        -\frac{1}{\phi^2} \log (p^q + (1-p)^q) = \frac{1}{\phi} \varphi_1 (q) \leq \tau_{\mu_p} (q) \leq \frac{1}{\phi-1} \varphi_1 (q) = -\log (p^q + (1-p)^q)\tc
    \end{equation*}
    where $\phi$ denotes the golden ratio.
    Reversing the inequalities yields the corresponding bounds for $q>1$.
    Of course, by considering larger $k$ we can obtain better bounds.
    For $p=1/2$, some bounds given by \cref{it:lq-inf-equiv} are shown in \cref{f:fib-bounds}.
\end{example}
\begin{figure}[t]
    \centering
    \begin{tikzpicture}[>=stealth,scale=1.8,font=\tiny]
    \begin{scope}[thick,gray,->]
        \draw (-0.1,0) -- (7.1,0);
        \draw (0,-0.7) -- (0,3);
    \end{scope}
    \foreach \y in {-1/2,1/2,1,3/2,2} \draw[gray] (-0.05,\y) node[black,left,font=\tiny]{$\y$}-- (0.05,\y);
    \foreach \x in {1,2,...,6} \draw[gray] (\x,-0.05) node[black,below,font=\tiny]{$\x$}-- (\x,0.05);

    \matrix [draw, fill=white, below right,scale=0.8,row sep=0.3em, column sep=0.3em, inner sep=0.3em, align=left] at (current bounding box.north west) {
        \draw[thick] (0,0) -- (0.7,0); & \node[anchor=west,inner sep=0pt] {$\varphi_k/(\lambda^k-1)$}; \\
        \draw[thick, dashed] (0,0) -- (0.7,0); & \node[anchor=west,inner sep=0pt] {$\varphi_k/\lambda^k$}; \\
    };

    \begin{scope}[thick,smooth, domain=0:7]
        \foreach\x in {3,5,7} {
            \draw[opacity={\x/7}] plot file {figures/fib_points_\x_lower.txt};
            \draw[dashed, opacity={\x/7}] plot file {figures/fib_points_\x_upper.txt};
        }
    \end{scope}
\end{tikzpicture}
    \caption{Upper and lower bounds on the $L^q$-spectrum of the frequency measure corresponding to the random Fibonacci substitution with $p=1/2$, for $k=3,5,7$.
        The darker shades correspond to higher values of $k$.
    }
    \label{f:fib-bounds}
\end{figure}

\begin{acknowledgements}
    The authors are grateful to Philipp Gohlke for his detailed comments on a draft version of this manuscript, which helped to remove some technical assumptions from \cref{it:multi-formalism}.
    They also thank Dan Rust and Tony Samuel for valuable input.
    AM thanks SFB1283 and the Universität Bielefeld for supporting a research visit during the summer of 2022, where some of the work on this project was undertaken.
    AM was supported by EPSRC DTP and the University of Birmingham.
    AR was supported by EPSRC Grant EP/V520123/1 and the Natural Sciences and Engineering Research Council of Canada.
    The authors thank the organisers of the Junior Ergodic Theory Meeting hosted at the ICMS in Edinburgh in March 2022, where this project began.
\end{acknowledgements}
\end{document}